\documentclass[final,leqno,onefignum,onetabnum]{amsart}

\newcommand{\mystrut}{\vrule height9.5pt depth1.5pt width0pt}
\newcommand{\tab}{\par\noindent\mystrut}
\newcommand{\tabb}{\tab\hskip1.5em}
\newcommand{\tabbb}{\tab\hskip3.5em}

\newcommand{\alg}[1]{\par\noindent\mystrut\ignorespaces\hbox to\textwidth{#1\hfill}}
\newcommand{\algt}[1]{\par\noindent\mystrut\hbox to\textwidth{\ignorespaces\hskip1.5em#1\hfill}}
\newcommand{\algtt}[1]{\par\noindent\mystrut\hbox to\textwidth{\ignorespaces\hskip3.5em#1\hfill}}
\newcommand{\algttt}[1]{\par\noindent\mystrut\hbox to\textwidth{\ignorespaces\hskip5.5em#1\hfill}}
\newcommand{\algtttt}[1]{\par\noindent\mystrut\hbox to\textwidth{\ignorespaces\hskip7.5em#1\hfill}}
\newcommand{\algttttt}[1]{\par\noindent\mystrut\hbox to\textwidth{\ignorespaces\hskip9.5em#1\hfill}}
\newcommand{\clg}[2]{\par\noindent\mystrut\hbox to\textwidth{\ignorespaces#1\hfill[#2]}}
\newcommand{\clgt}[2]{\par\noindent\mystrut\hbox to\textwidth{\ignorespaces\hskip1.5em#1\hfill[#2]}}
\newcommand{\clgtt}[2]{\par\noindent\mystrut\hbox to\textwidth{\ignorespaces\hskip3.5em#1\hfill[#2]}}
\newcommand{\clgttt}[2]{\par\noindent\mystrut\hbox to\textwidth{\ignorespaces\hskip5.5em#1\hfill[#2]}}
\newcommand{\clgtttt}[2]{\par\noindent\mystrut\hbox to\textwidth{\ignorespaces\hskip7.5em#1\hfill[#2]}}
\newcommand{\clgttttt}[2]{\par\noindent\mystrut\hbox to\textwidth{\ignorespaces\hskip9.5em#1\hfill[#2]}}

\newcommand{\FOR}{\textbf{for}\hskip2pt}
\newcommand{\IF}{\textbf{if}\hskip2pt}

\newcommand{\END}{\textbf{end}}

\newcommand{\ELSE}{\hskip2pt\textbf{else}\hskip2pt}

\newcommand{\ELSEIF}{\textbf{else if}}

\usepackage{array,latexsym,pst-node,amsmath,amssymb}
\usepackage{amsfonts,cite}
\usepackage{graphicx}
\usepackage{bbm}
\usepackage{multirow}

\usepackage[norelsize]{algorithm2e}  

\newcolumntype{x}[1]{%
&gt;{\raggedleft\arraybackslash}p{#1}}%

\title{On solving large-scale limited-memory quasi-Newton equations}

\author{Jennifer B. Erway}
\email{erwayjb@wfu.edu}
\address{Department of Mathematics, PO Box 7388,
    Wake Forest University, Winston-Salem, NC 27109}
\author{Roummel F. Marcia}
\email{rmarcia@ucmerced.edu}
\address{School of Natural Sciences, University of California,
    Merced, 5200 N. Lake Road, Merced, CA 95343}
\thanks{Research supported in part by NSF grants CMMI-1334042 and CMMI-1333326.}

\begin{document}

\newcommand{\subject}{\mathop{\mathrm{subject\ to}}}  
\newcommand{\minimize}[1]{{\displaystyle\minim_{#1}}}

\newcommand{\minim}{\mathop{\mathrm{minimize}}}
\newcommand{\minimizearg}[2]{{\displaystyle\minim_{#1 \in \Re^{#2}}}}
\newcommand{\maximizearg}[2]{{\displaystyle\maxim_{#1 \in \Re^{#2}}}}
\newcommand{\minimizetwoarg}[4]{{\displaystyle\minim_{#1 \in \Re^{#2},%
                                                  #3 \in \Re^{#4}}}}
\newcommand{\maximize}[1]{{\displaystyle\maxim_{#1}}}
\newcommand{\maxim}{\mathop{\mathrm{maximize}}}
\renewcommand{\maximize}[1]{{\displaystyle\maxim_{#1}}}
\newcommand{\MATLAB}{{\small MATLAB}}
\newcommand{\SMW}{{\small SMW}}
\newcommand{\QR}{{\small QR}}
\newcommand{\BFGS}{{\small BFGS}}
\newcommand{\DFP}{{\small DFP}}
\newcommand{\SR}{{\small SR1}}
\newcommand{\LSR}{{\small L-SR1}}
\newcommand{\LBFGS}{{\small L-BFGS}}
\newcommand{\mgap}{\;\;}
\newcommand{\bgap}{\;\;\;}
\newcommand{\defined}{\mathop{\,{\scriptstyle\stackrel{\triangle}{=}}}\,}

\makeatletter
\newcommand{\infimum}[1]{{\displaystyle\infim_{#1}}}
\newcommand{\infim}{\mathop{\operator@font{inf}}}
\newcommand{\todo}[1]{{\textcolor{red}{[#1]}}}
\newcommand{\rfm}[1]{{\textcolor{blue}{#1}}}
\newcommand\XOR{\mathbin{\char`\^}}

\newcommand{\algbox}{%
  \vbox{\hrule height0.3pt\hbox{\vrule height2pt%
      width0.3pt\hskip\textwidth\hskip-0.6pt\vrule width0.3pt}\hrule height0.3pt}}

\newcommand{\AlgBegin}{\vspace{\smallskipamount}\noindent\algbox\par\vspace{\smallskipamount}\par\noindent}
\newcommand{\AlgEnd}{\par\noindent\algbox\vspace{\smallskipamount}}

\newcounter{pseudocode}[section]
\def\thepseudocode{\thesection.\arabic{pseudocode}}
\newenvironment{pseudocode}[2]%
        {%
        \refstepcounter{pseudocode}%
          \AlgBegin %
               {{\bfseries Algorithm \thepseudocode.}\rule[-1.25pt]{0pt}{10pt}#1}%
        #2}%
           {\AlgEnd}

\newcounter{Pseudocode}[section]
\def\thePseudocode{\thesection.\arabic{Pseudocode}}
\newenvironment{Pseudocode}[2]%
        {%
        \refstepcounter{Pseudocode}%
          \AlgBegin %
               {{\bfseries #1.}\rule[-1.25pt]{0pt}{10pt}}%
        #2}%
           {\AlgEnd}

\def\tnu{\tilde{\nu}}
\makeatother

\maketitle

\begin{abstract}
We consider the problem of solving linear systems of equations arising with 
limited-memory members of the restricted Broyden class of updates
and the symmetric rank-one (SR1) update.  
In this paper, we propose a new
approach based on a practical implementation of the compact
representation for the inverse of these limited-memory matrices.
Numerical results suggest
that the proposed method compares favorably in speed and accuracy to
other algorithms and is competitive with several update-specific
methods available to only a few members of the Broyden class of updates.
Using the proposed approach has an additional benefit: The condition number of
the system matrix can be computed efficiently.  
\end{abstract}


\pagestyle{myheadings}
\thispagestyle{plain}
\markboth{J. B. ERWAY AND R. F. MARCIA}{ON SOLVING LARGE-SCALE LIMITED-MEMORY QUASI-NEWTON EQUATIONS}

\section{Introduction}
We consider linear systems of the following form:
\begin{equation} \label{eqn-basic}
	B_{k+1} r = z,
\end{equation}
where $r,z\in\Re^n$ and $B_{k+1} \in \Re^{n \times n}$ is a limited-memory
quasi-Newton matrix obtained from applying $k+1$ Broyden class updates to
an initial matrix $B_0$.  We assume $n$ is large, and thus, explicitly
forming and storing $B_{k+1}$ is impractical or impossible; moreover, in
this setting, we assume limited-memory quasi-Newton matrices so that only
the most recently-computed $k$ updates are stored and used to update $B_0$.
In practice, the value of $k$ is small, i.e., less than 10 (see e.g.,
~\cite{ByrNS94}), making $k\ll n$.  Problems such as (\ref{eqn-basic})
arise in quasi-Newton line-search and trust-region methods for large-scale
optimization (see, e.g., ~\cite{ConGT00a,DenM77,GrNaS09,NocW06}), as well
as in preconditioning iterative solvers (see, e.g., ~\cite{MorN00,OLeY94}).
In this paper, we propose a compact formulation of $B_{k+1}^{-1}$ that can
be used to efficiently solve (\ref{eqn-basic}).

\bigskip

Traditional quasi-Newton methods for minimizing a continuously
differentiable function $f:\Re^n\rightarrow \Re$ generate a
sequence of iterates $\{x_k\}$ such that $f$ is strictly decreasing on this
sequence.  Moreover, at each iteration, the most-recently
computed iterate $x_{k+1}$ is used to update the quasi-Newton matrix by defining
a new \emph{quasi-Newton pair} $(s_k, y_k)$ given by
$$s_k\defined x_{k+1}-x_k \quad \text{and} \quad y_k\defined \nabla f(x_{k+1})-\nabla f(x_k).$$
In the case of the Broyden class updates, $B_{k+1}$
is updated as follows:
\begin{equation}\label{eqn-broyden}
 B_{k+1}  =B_k   -  \frac{1}{s_k ^TB_k  s_k }B_k  s_k s_k ^TB_k  +  
  	\frac{1}{y_k^Ts_k }y_ky_k^T +
         \phi (s_k^T B_k  s_k )w_kw_k^T,
\end{equation}
with
$$
	w_k = \frac{y_k}{y_k^Ts_k} - \frac{B_k  s_k}{s_k^TB_k  s_k},
$$
where $\phi\in\Re$.  The most well-known quasi-Newton update is the
Broyden-Fletcher-Goldfarb-Shanno (\BFGS) update, which is obtained by
setting $\phi=0$.  While this is the most widely-used update, research has
suggested that other values of $\phi$ may lead to faster
convergence~\cite{preconvex,byrd1992behavior,liu2007statistical}.  
Because of this interest in other values of $\phi$, the entire Broyden
class (i.e., $\phi\in\Re$) of quasi-Newton methods has been
generalized to solve minimization problems over Riemannian
manifolds~\cite{huang2015broyden}.  For these reasons, in this paper,
we consider solving (\ref{eqn-basic}) for other values of $\phi$ other
than $\phi=0$ (the BFGS update).  In particular, we consider the
\emph{restricted} Broyden class of updates, where $\phi \in [0, 1]$.  Under
mild assumptions, quasi-Newton matrices in this class are positive
definite, a property that is important for computing search directions in
optimization.  We also consider the special case of the symmetric rank-one
(SR1) update, an update that has been of recent interest in
trust-region methods for unconstrained optimization
\cite{BruEM15,BruBEMY16}.  (For more details on
limited-memory quasi-Newton matrices and the Broyden class of matrices see,
e.g.,~\cite{ByrNS94,GrNaS09, NocW06}.)

\bigskip

Solving systems of the form (\ref{eqn-basic}) can be done efficiently in
the case of the \BFGS{} update ($\phi=0$) using the well-known two-loop
recursion~\cite{Noc80} (see \ref{app:two-loop} for details).  However,
there is no known corresponding recursion method for other updates of the
Broyden class (see \cite[p.190]{NocW06}); in fact, for this reason many
researchers prefer the \BFGS{} update.  In this paper we present various
methods for solving linear systems of equations with limited-memory
quasi-Newton matrices.  The main contribution of this paper is a new
approach by formulating a compact representation for the inverse of any
member of the restricted Broyden class.  This representation
allows us to efficiently solve linear
systems of the form (\ref{eqn-basic}).   
The compact formulation for the
inverse of these matrices is based on ideas found in~\cite{ErwM15}, where a
compact formulation for the restricted Broyden class of matrices is
presented.  
An additional benefit of our proposed approach is the ability
to calculate the eigenvalues of the limited-memory matrix, and hence, the
condition number of the linear system.

This paper is organized in seven main sections.  In Section 2, we present
current methods for solving linear systems with a limited-memory
quasi-Newton matrix.  In Section 3, we review the compact formulation for
matrices obtained using the restricted Broyden class of updates.  The
compact formulation for the inverse of these matrices is presented in
Section 4.  In Section 5, we provide a practical implementation for solving
systems of the form (\ref{eqn-basic}) when the system matrix is either
  a member of the restricted Broyden class or an \SR{} matrix.  In
addition, we discuss how to compute the condition number of the linear
system being solved and how to obtain additional computational savings when
a new quasi-Newton pair is computed.  We present numerical experiments in
Section 6 that demonstrate the competitiveness of our proposed approach.
Finally, concluding remarks are in Section 7.

\section{Current methods}

One approach to solve linear equations with members of the Broyden class is
to use the Sherman-Morrison-Woodbury (\SMW) formula to update the inverse
after each rank-one change (see, e.g., ~\cite{GolV96}).  Thus, to compute
the inverse after a rank-two update, the \SMW{} formula 
may be applied 
twice.  (For algorithmic details on this approach, see \ref{app:SMW-recursion}.)
Alternatively, linear solves with a general member of the Broyden class can
be performed using the following recursion formula found in~\cite{DenM77}
for $H_{k+1}$, the inverse of $B_{k+1}$.  Specifically, the inverse of
  $B_{k+1}$ in \eqref{eqn-broyden} is given by the following rank-two
  update:
\begin{equation}\label{eqn-1paraminv}
	H_{k+1}  = H_k  + \frac{1}{s_k^Ty_k}s_ks_k^T 
				- \frac{1}{y_k^TH_k y_k}H_k y_ky_k^TH_k 
				+ \Phi_k (y_k^TH_k y_k)v_kv_k^T,
\end{equation}
where
\begin{equation}\label{eqn-vPHI}
	v_k = \frac{s_k}{y_k^Ts_k} - \frac{H_k y_k}{y_k^TH_k y_k},
	\
	\text{and} 
	\ \
	\Phi_k 
	= 
	\frac{(1-\phi)(y_k^Ts_k)^2}{(1-\phi)(y_k^Ts_k)^2 + \phi(y_k^TH_k y_k)(s_k^TB_k s_k)}, \ \ 
\end{equation}
and $H_k\defined B_k^{-1}$.  (Algorithmic details associated with this
approach can be found in \ref{app:Hk+1}.)  

The method proposed in this
paper makes use of the so-called compact formulation of quasi-Newton
matrices that not only appears to be competitive in terms of speed but also
yields additional information about the system matrix.

\section{Compact representation}
In this section, we briefly review the compact representation of 
members of the Broyden class of quasi-Newton matrices.

Given $k+1$ update pairs of the form $\{(s_i,y_i)\}$,
the \emph{compact formulation} of $B_{k+1}$ is given by
\begin{equation}\label{eqn-compact}
	B_{k+1}  = B_0  + \Psi_k M_k \Psi_k^T,
\end{equation}
where $M_k$ is a square matrix
and $B_0\in\Re^{n\times n}$ is an initial matrix, which is often taken to
be a scalar multiple of the identity.  (Note that the size of $\Psi$ is not
specified.)  The compact formulation of members of the Broyden class of
updates are defined in terms of
\begin{eqnarray*}
	S_k &=& [ \ s_0 \ \ s_1 \ \ s_2 \ \ \cdots \ \ s_{k} \ ] \ \in \ \Re^{n \times (k+1)}, \\
	Y_k &=& [ \ y_0 \ \ y_1 \ \ y_2 \ \ \cdots \ \ y_{k} \ ] \ \in \ \Re^{n \times (k+1)},
\end{eqnarray*}
and the following decomposition of $S_k^TY_k \in\Re^{(k+1) \times (k+1)}$:
$$
	S_k^TY_k =   L_k + D_k + R_k,
$$
where $L_k$ is strictly lower triangular, $D_k$ is diagonal, and $R_k$ is
strictly upper triangular.  

Assuming all updates are well defined,
Erway and Marcia~\cite{ErwM15} derive the compact formulation
for all members of the Broyden class for $\phi\in [0,1]$.  In particular,
for any $\phi\in [0,1]$, the compact formulation for $B_{k+1}$ is given
by (\ref{eqn-compact}) 
where
\begin{equation}\label{eqn-alt-form-convex}
	\Psi_k = \begin{bmatrix} B_0S_k & Y_k \end{bmatrix} \quad \text{and} \quad
	M_k  \ = \
	-\begin{bmatrix}
	S_k^TB_0S_k - \phi \Lambda_k & \ \ (L_k-\phi \Lambda_k) \\
	(L_k-\phi \Lambda_k)^T & -(D_k + \phi \Lambda_k) 
	\end{bmatrix}^{-1},
\end{equation}
where
\begin{equation}\label{eqn-Lambda}
	\Lambda_k = \underset{0 \le i \le k}{\text{diag}} \ (\lambda_i), \quad \text{where }
	\lambda_i = 
	\frac{1}
	{\displaystyle
	-\frac{1-\phi}{s_i^TB_is_i} 
	-\frac{\phi}{s_i^Ty_i}}
	\ \ \text{for $0 \le i \le k$.}
\end{equation}
When $\phi=0$, this representation simplifies to the \BFGS{} compact
representation found in~\cite{ByrNS94}.  In the limited-memory case, $k\ll
n$, and thus, $M_k$ is easily inverted.

\bigskip

\section{Inverses}

In this section, we present the compact formulation for the inverse
of any member of the restricted Broyden class of matrices and for \SR{}
matrices.  We then demonstrate how this compact formulation can be 
used to solve linear systems with these limited-memory matrices.
Finally, we discuss computing the condition number of the linear system and
potential computational savings when a new quasi-Newton pair is computed.

The main contribution of this paper is formulating a compact
  representation for the inverse for any member of the restricted Broyden
class.  It should be noted that the compact formulation for the inverse of
a \BFGS{} matrix is already known (see ~\cite{ByrNS94}).  We now derive a
general expression for the compact representation of any member of the
restricted Broyden class.  To do this, we apply the
Sherman-Morrison-Woodbury formula (see, e.g.,~\cite{GolV96}) to the compact
representation of $B_{k+1}$ found in (\ref{eqn-compact}):
$$
	B_{k+1}^{-1}= 
	B_0^{-1} + B_0^{-1} \Psi_k \left (-M_k^{-1} - \Psi_k^TB_0^{-1} \Psi_k \right )^{-1} \Psi_k^T B_0^{-1}.
$$
For quasi-Newton matrices it is conventional to let $H_i$ denote the inverse of $B_i$ for each $i$;
using this notation, the inverse of $B_{k+1}^{-1}$ is given by
\begin{equation}\label{eqn-H1}
	H_{k+1}= 
	H_0 + H_0 \Psi_k \left (-M_k^{-1} - \Psi_k^TH_0 \Psi_k \right )^{-1} \Psi_k^T H_0.
\end{equation}
Since
$ H_0 \Psi_k 
	=
	H_0 [ B_0 S_k \,\,\, Y_k ]  
	=
	[ S_k \,\,\, H_0 Y_k  ],
$ (\ref{eqn-H1}) can be written as
\begin{equation}\label{eqn-Hk+1compact}
	H_{k+1} =
	H_0 +
	[ \ S_k \ \ \ H_0 Y_k \ ]
	(-M_k^{-1} - \Psi_k^TH_0\Psi_k)^{-1} 
	\begin{bmatrix}
		S_k^T \\
		Y_k^TH_0
	\end{bmatrix}.
\end{equation}
Finally, (\ref{eqn-Hk+1compact}) can be simplified
using
$$
	\Psi_k^T H_0 \Psi_k =
	\begin{bmatrix}
		S_k^TB_0S_k &  S_k^TY_k \\
		Y_k^TS_k & Y_k^TH_0Y_k
	\end{bmatrix}
$$
to obtain
\begin{eqnarray*}
	&& \hspace{-1cm} (-M_k^{-1} - \Psi_k^TH_0\Psi_k)^{-1} \\
	\ \ \ &=&
	\left (
	\begin{bmatrix}
	S_k^TB_0S_k - \phi \Lambda_k & \ \ (L_k-\phi \Lambda_k) \\
	(L_k-\phi \Lambda_k)^T & -(D_k + \phi \Lambda_k) 
	\end{bmatrix}
	-
	\begin{bmatrix}
		S_k^TB_0S_k &  S_k^TY_k \\
		Y_k^TS_k & Y_k^TH_0Y_k
	\end{bmatrix}
	\right )^{-1}
	\\
	&=&
	\begin{bmatrix}
	-\phi \Lambda_k & L_k - \phi \Lambda_k - S_k^TY_k  \\
	L_k^T - \phi \Lambda_k - Y_k^TS_k & -D_k - \phi\Lambda_k - Y_k^TH_0Y_k 
	\end{bmatrix}^{-1}\\[.2cm]
	&=&
	\begin{bmatrix}
	-\phi \Lambda_k & -R_k - D_k - \phi \Lambda_k \\
	-R_k^T - D_k - \phi \Lambda_k  & -D_k - \phi\Lambda_k - Y_k^TH_0Y_k 
	\end{bmatrix}^{-1}.
\end{eqnarray*}
In other words, the compact formulation for the inverse of a member
of the restricted Broyden class is given by
\begin{equation}\label{eqn-H2}
	H_{k+1}= 
	H_0 + \tilde{\Psi}_k \tilde{M}_k \tilde{\Psi}_k^T \quad
\text{where} \quad \tilde{\Psi}_k \defined [S_k \,\,\, H_0Y_k], 
\end{equation}
and 
\begin{equation}\label{eqn-mtilde}
\tilde{M}_k \defined 
	\begin{bmatrix}
	-\phi \Lambda_k & -R_k - D_k - \phi \Lambda_k \\
	-R_k^T - D_k - \phi \Lambda_k  & -D_k - \phi\Lambda_k - Y_k^TH_0Y_k 
	\end{bmatrix}^{-1}.
\end{equation}

\bigskip

In the following section, we discuss a practical method to compute
$\tilde{M}_k$.  In \ref{app:updates}, we explore remarkable
relationships between various updates.

\section{Recursive formulation and practical implementation}

In this section, we present a practical method to compute
$\tilde{M}_k$ in (\ref{eqn-H2}) given by (\ref{eqn-mtilde}).  
We begin by providing an alternative expression for $\tilde{M}_0$
that will allow us to define a recursion method to compute $\tilde{M}_k$.

\bigskip

\noindent
\textbf{Lemma 1.}
\textsl{
Suppose $H_{1}=H_0+\tilde{\Psi}_0\tilde{M}_0\tilde{\Psi}_0^T$,
is the inverse of
a member of the restricted Broyden class of updates after performing
one update.  Then,
\begin{equation}\label{eqn-M0}
\tilde{M}_0 = 
        \begin{bmatrix}
              \tilde{\alpha}_0 & \tilde{\beta}_0 \\ \tilde{\beta}_0 & \tilde{\delta}_0
        \end{bmatrix} , 
\end{equation}
where
\begin{equation}\label{eqn-base0}
	\tilde{\alpha}_0 = \frac{1}{s_0^Ty_0}+ \Phi_0\frac{y_0^TH_0y_0}{(s_0^Ty_0)^2}, \quad
	\tilde{\beta}_0 = -\frac{\Phi_0}{y_0^Ts_0}, \quad
	\tilde{\delta}_0 = -\frac{1-\Phi_0}{y_0^TH_0y_0},
\end{equation}
and $\Phi_0$ is given in \eqref{eqn-vPHI}.
}

\bigskip

\noindent \textbf{Proof.}
Expanding \eqref{eqn-1paraminv}, yields
\begin{eqnarray*}
        H_1&=& 
        H_0+
        \left (
        		\frac{1}{s_0^Ty_0} + \Phi_0\frac{y_0^TH_0y_0}{(s_0^Ty_0)^2}
        \right ) 
        s_0s_0^T
        - \frac{\Phi_0}{y_0^Ts_0}H_0y_0s_0^T \\
        &&
        \qquad - \frac{\Phi_0}{y_0^Ts_0}s_0y_0^TH_0 
        - \frac{1 - \Phi_0}{y_0^TH_0y_0} H_0y_0y_0^TH_0,
\end{eqnarray*}
which simplifies to
\begin{equation}\label{eqn-H0+1}
H_1        = H_0 + 
                \begin{bmatrix}
                s_0 & H_0y_0
        \end{bmatrix}
        \begin{bmatrix}
              \tilde{\alpha}_0 & \tilde{\beta}_0 \\ \tilde{\beta}_0 & \tilde{\delta}_0
        \end{bmatrix}
        \begin{bmatrix}
                s_0^T \\
                y_0^TH_0
        \end{bmatrix},
\end{equation}
where
$\tilde{\alpha}_0$, $\tilde{\beta}_0$, and $\tilde{\delta}_0$
are defined as in (\ref{eqn-base0})
and $\Phi_0$ is given by (\ref{eqn-vPHI}).
Note the \eqref{eqn-H0+1} is of the form
$H_{1}=H_0+\tilde{\Psi}_0\tilde{M}_0\tilde{\Psi}_0^T$.

We now show that $\tilde{M}_0$ defined by 
(\ref{eqn-M0}) and (\ref{eqn-base0}) is equivalent to
(\ref{eqn-mtilde}) with $k=0$.
To see this, we simplify the entries of (\ref{eqn-M0}).
First, define
$$
	\Delta_0 \defined (1-\phi)(y_0^Ts_0)^2 + \phi (y_0^TH_0y_0)(s_0^TB_0s_0).
$$
Note that $\Delta_0 \ne 0$ since $H_0$ and $B_0$ are positive definite
and $\phi \in [0, 1]$.  Then,
$\Phi_0 = (1-\phi)(s_0^Ty_0)^2/\Delta_0,$
and thus,
$$
	\tilde{\alpha}_0 = \frac{1}{s_0^Ty_0} + \frac{(1-\phi)y_0^TH_0y_0}{\Delta_0}, \ \
	\tilde{\beta}_0 = - \frac{(1-\phi)(s_0^Ty_0)}{\Delta_0}, \ \ \text{and} \ \
	\tilde{\delta}_0 = -\frac{\phi (s_0^TB_0s_0)}{\Delta_0}.
$$
Consider the inverse of (\ref{eqn-M0}), which is given by
\begin{equation}\label{eqn-base0inv}
       \begin{bmatrix}
             \tilde{\alpha}_0 & \tilde{\beta}_0 \\ \tilde{\beta}_0 & \tilde{\delta}_0
       \end{bmatrix}^{-1}
       =
	\frac{1}{\tilde{\alpha}_0\tilde{\delta}_0-\tilde{\beta}_0^2}
       \begin{bmatrix}
             \tilde{\delta}_0 & -\tilde{\beta}_0 \\ -\tilde{\beta}_0 & \tilde{\alpha}_0
       \end{bmatrix}.
\end{equation}
We now simplify the entries on the right side of (\ref{eqn-base0inv}).
The determinant of (\ref{eqn-M0}) is given by
\begin{eqnarray}
	\!\!\!\!\!\!
	\tilde{\alpha}_0\tilde{\delta}_0 - \tilde{\beta}_0^2 
	\! = \!
	\frac{-(1-\Phi_0)}{(s_0^Ty_0)(y_0^TH_0y_0)}
	-\frac{\Phi_0(1-\Phi_0)}{(s_0^Ty_0)^2}
	-\frac{\Phi_0^2}{(s_0^Ty_0)^2}
	\! = \!
	 \frac{1}{\Delta_0} \left (
	\frac{s_0^TB_0s_0}{\lambda_0}
	\right )\!, \label{eqn-determtilde}
\end{eqnarray}
where $\lambda_0$ in \eqref{eqn-determtilde} is defined in \eqref{eqn-Lambda}.
Thus, the first entry of (\ref{eqn-base0inv}) is given by
\begin{equation}\label{eqn-deltatilde}
	\frac{\tilde{\delta}_0}{\tilde{\alpha}_0\tilde{\delta}_0 - \tilde{\beta}_0^2} 
	=
	-\frac{\phi (s_0^TB_0s_0)}{\Delta_0}
	\frac{\Delta_0\lambda_0}{s_0^TB_0s_0} =
-\phi \lambda_0.\end{equation}
The off-diagonal elements of the right-hand side of (\ref{eqn-base0inv})
simplify as follows:
\begin{eqnarray}
	\frac{\tilde{\beta}_0}{\tilde{\alpha}_0\tilde{\delta}_0 - \tilde{\beta}_0^2}
	=
	- \left (\frac{(1-\phi)(s_0^Ty_0)}{\Delta_0} \right )
	\left (
	{\Delta_0} 
	\left (
	\frac{\lambda_0}{s_0^TB_0s_0}
	\right )
	\right )
	=
	s_0^Ty_0 + \phi \lambda_0. \label{eqn-betaguys} 
\end{eqnarray}
Finally, the last entry of (\ref{eqn-base0inv}) can be simplified as
follows:
\begin{eqnarray}
	\nonumber
	\frac{\tilde{\alpha}_0}{\tilde{\alpha}_0\tilde{\delta}_0 - \tilde{\beta}_0^2} 
	&=&
	\frac{1}{s_0^Ty_0} \frac{\Delta_0 \lambda_0}{s_0^TB_0s_0} +
	(1-\phi)y_0^TH_0y_0\frac{\lambda_0}{s_0^TB_0s_0}\\
	&=& -s_0^Ty_0 - \phi\lambda_0 - y_0^TH_0y_0 . \label{eqn-alphaguy}
\end{eqnarray}
(For details in the calculations of \eqref{eqn-determtilde},  \eqref{eqn-betaguys}, \eqref{eqn-alphaguy}, see 
 \ref{app:Lemma1}.)
Thus, (\ref{eqn-base0inv}) together with
(\ref{eqn-deltatilde}),
(\ref{eqn-betaguys}), and (\ref{eqn-alphaguy})
gives
$$
        \begin{bmatrix}
              \tilde{\alpha}_0 & \tilde{\beta}_0 \\ \tilde{\beta}_0 & \tilde{\delta}_0
        \end{bmatrix}^{-1} = 
        \begin{bmatrix}
        		-\phi \lambda_0 & -s_0^Ty_0 - \phi \lambda_0  \\
		-s_0^Ty_0 - \phi \lambda_0 & -s_0^Ty_0 - \phi\lambda_0 - y_0^TH_0y_0
       \end{bmatrix},
$$
showing that
$\tilde{M}_0$ defined by 
(\ref{eqn-M0}) and (\ref{eqn-base0}) is equivalent to
(\ref{eqn-mtilde}) with $k=0$.
$\square$

\bigskip

Together with Lemma 1, the following theorem shows how $\tilde{M}_k$ is
related to $\tilde{M}_{k-1}$; from this, we present an algorithm to compute
$\tilde{M}_k$ using recursion.  This computation avoids explicitly forming
$B_k$, which is used in the definition of $\Lambda_k$ and $\Phi_k$.

\bigskip

\noindent
\textbf{Theorem 1.}
\textsl{Suppose $H_{j+1}=H_0+\tilde{\Psi}_j\tilde{M}_j\tilde{\Psi}_j^T$
as in (\ref{eqn-H2}) for all $j \in \{0, \dots, k\}$, is the inverse of
a member of the Broyden class of updates for a fixed $\phi\in [0,1]$.  
If $\tilde{M}_0$ is defined by (\ref{eqn-M0}), then
for all $j\in\{1,\ldots k\}$,
$\tilde{M}_j$ satisfies the recursion relation
\begin{equation}\label{eqn-th1}
\tilde{M}_j=
	\Pi_j^T \!\!
	\begin{bmatrix}
		\tilde{M}_{j-1}
			+\tilde{\delta}_j\tilde{u}_j\tilde{u}_j^T & \tilde{\beta}_j\tilde{u}_j 
			& \tilde{\delta}_j \tilde{u}_j \\
		\tilde{\beta}_j\tilde{u}_j^T & \tilde{\alpha}_j & \tilde{\beta}_j \\
		\tilde{\delta}_j\tilde{u}_j^T &\tilde{\beta}_j & \tilde{\delta}_j
	\end{bmatrix}
	\! \Pi_j,
	\ \text{where} \
	\Pi_j \! = \!\!
	\begin{bmatrix}
		I_{j} & 0 & 0 & 0 \\
		0 & 0 & 1 & 0 \\
		0 & I_{j} & 0 & 0 \\
		0 & 0 & 0 & 1
	\end{bmatrix} 
\end{equation}
is such that
$
	\begin{bmatrix}
		\tilde{\Psi}_{j-1} & s_j & H_0y_j 
	\end{bmatrix}
	\Pi_j
	= \tilde{\Psi}_j,
$
and
\begin{equation*}
	\tilde{\alpha}_j = \frac{1}{s_j^Ty_j}+ \Phi_j\frac{y_j^TH_jy_j}{(s_j^Ty_j)^2}, \quad
	\tilde{\beta}_j = -\frac{\Phi_j}{y_j^Ts_j}, \quad
	\tilde{\delta}_j = -\frac{1-\Phi_j}{y_j^TH_jy_j},
\end{equation*}
and $\tilde{u}_j = \tilde{M}_{j-1}\tilde{\Psi}_{j-1}^Ty_j$.}

\bigskip

\noindent \textbf{Proof.}
This proof is by induction on $j$.  For the base case, we
show that 
equation (\ref{eqn-th1}) holds for $j=1$.
Setting $k=1$ in
\eqref{eqn-1paraminv} yields
\begin{equation}\label{eqn-H_2}
	H_{2} = H_1 + 
	 \begin{bmatrix}
                s_1 & H_1y_1
        \end{bmatrix}
        \begin{bmatrix}
              \tilde{\alpha}_1 & \tilde{\beta}_1 \\ \tilde{\beta}_1 & \tilde{\delta}_1
        \end{bmatrix}
        \begin{bmatrix}
                s_1^T \\
                y_1^TH_1
        \end{bmatrix},
\end{equation}
where
\begin{equation*}
	\tilde{\alpha}_1 = \frac{1}{s_1^Ty_1}+ \Phi_1\frac{y_1^TH_1y_1}{(s_1^Ty_1)^2}, \quad
	\tilde{\beta}_1 = -\frac{\Phi_1}{y_1^Ts_1}, \quad \text{and} \quad
	\tilde{\delta}_1 = -\frac{1-\Phi_1}{y_1^TH_1y_1}.
\end{equation*}
By Lemma 1, (\ref{eqn-H_2}) can be written as
\begin{equation}\label{eqn-thrm1base}
H_2  =
	H_0 + \tilde{\Psi}_{0}\tilde{M}_{0}\tilde{\Psi}_{0}^T + 
	 \begin{bmatrix}
                s_1 & H_1y_1
        \end{bmatrix}
        \begin{bmatrix}
              \tilde{\alpha}_1 & \tilde{\beta}_1 \\ \tilde{\beta}_1 & \tilde{\delta}_1
        \end{bmatrix}
        \begin{bmatrix}
                s_1^T \\
                y_1^TH_1
        \end{bmatrix},
\end{equation}
where $\tilde{\Psi}_0 = [S_0 \,\,\, H_0Y_0] $ and $\tilde{M}_0$ is given by (\ref{eqn-M0}).
Letting $\tilde{u}_1 = \tilde{M}_{0}\tilde{\Psi}_{0}^Ty_1$, we have that
\begin{equation}\label{eqn-H1y1}
	H_1y_1 = \left ( H_0 +  \tilde{\Psi}_{0}\tilde{M}_{0}\tilde{\Psi}_{0}^T  \right ) y_1 = H_0y_1 
			+ \tilde{\Psi}_{0}\tilde{u}_{1}.
\end{equation}
Substituting (\ref{eqn-H1y1}) into the last quantity on the
right side of (\ref{eqn-thrm1base}) yields
\begin{eqnarray*}
	\begin{bmatrix}
                s_1 & H_1y_1
        \end{bmatrix}
        \begin{bmatrix}
              \tilde{\alpha}_1 & \tilde{\beta}_1 \\ \tilde{\beta}_1 & \tilde{\delta}_1
        \end{bmatrix}
        \begin{bmatrix}
                s_1^T \\
                y_1^TH_1
        \end{bmatrix}
 	&=&
	\begin{bmatrix}
                s_1 & H_0y_1 +  \tilde{\Psi}_{0}\tilde{u}_{1}
        \end{bmatrix}
        \begin{bmatrix}
              \tilde{\alpha}_1 & \tilde{\beta}_1 \\ \tilde{\beta}_1 & \tilde{\delta}_1
        \end{bmatrix}
        \begin{bmatrix}
                s_1^T \\
                y_1^TH_0 + \tilde{u}_{1}^T \tilde{\Psi}_{0}^T
        \end{bmatrix}\\
        &=&
        \begin{bmatrix}
        		\tilde{\Psi}_{0} & s_1 & H_0y_1
        \end{bmatrix}
	\begin{bmatrix}
		\tilde{\delta}_1\tilde{u}_1\tilde{u}_1^T & \tilde{\beta}_1\tilde{u}_1 
			& \tilde{\delta}_1 \tilde{u}_1 \\
		\tilde{\beta}_1\tilde{u}_1^T & \tilde{\alpha}_1 & \tilde{\beta}_1 \\
		\tilde{\delta}_1\tilde{u}_1^T &\tilde{\beta}_1 & \tilde{\delta}_1
	\end{bmatrix}
        \begin{bmatrix}
        		\tilde{\Psi}_{0}^T \\
		s_1^T \\
		y_1^TH_0
        \end{bmatrix}.
\end{eqnarray*}
Letting $\Pi_1$ be defined as in \eqref{eqn-th1} with $j=1$, gives that
$
	\begin{bmatrix}
		\tilde{\Psi}_{0} & s_1 & H_0y_1 
	\end{bmatrix}
	\Pi_1
	= \tilde{\Psi}_1,
$ and thus,
$$
	H_2 = H_0 +
	\tilde{\Psi}_1
	\Pi_1^T
        \bar{M}_1
	\Pi_1 \tilde{\Psi}_1^T,
$$
where $$\bar{M}_1\defined
	\begin{bmatrix}
		\tilde{M}_{0}
			+\tilde{\delta}_1\tilde{u}_1\tilde{u}_1^T & \tilde{\beta}_1\tilde{u}_1 
			& \tilde{\delta}_1 \tilde{u}_1 \\
		\tilde{\beta}_1\tilde{u}_1^T & \tilde{\alpha}_1 & \tilde{\beta}_1 \\
		\tilde{\delta}_1\tilde{u}_1^T &\tilde{\beta}_1 & \tilde{\delta}_1
	\end{bmatrix}.
$$
It remains to show that $\tilde{M}_1=\Pi_1^T\bar{M}_1\Pi_1$; however, for
simplicity, we instead show $\tilde{M}_1^{-1} = \Pi_1^T\bar{M}_1^{-1}\Pi$.

Notice that $\bar{M}_1$ can be written
as the product of three matrices:
\begin{eqnarray*}
\bar{M}_1		&=&
		\begin{bmatrix}
			I_2 & 0 & \tilde{u}_1\\
			0 & 1 & 0 \\
			0 & 0 & 1
		\end{bmatrix}
		\begin{bmatrix}
			\tilde{M}_{0} & 0   \\
			0 & \tilde{N}_0 
		\end{bmatrix}
		\begin{bmatrix}
			I_2 & 0 & 0 \\
			0 & 1 & 0 \\
			\tilde{u}_1^T & 0 & 1
		\end{bmatrix},
		\ \text{where } 
		\tilde{N}_0=\begin{bmatrix}
		\tilde{\alpha}_1 & \tilde{\beta}_1 \\
		 \tilde{\beta}_1&  \tilde{\delta}_1
		\end{bmatrix}.
\end{eqnarray*}
Then, $\bar{M}_1^{-1}$ is given by
\begin{eqnarray}\label{eqn-biginv}
\bar{M}_1^{-1}		&=&
		\begin{bmatrix}
			I & 0 & 0 \\
			0 & 1 & 0 \\
			-\tilde{u}_1^T & 0 & 1
		\end{bmatrix}
		\begin{bmatrix}
			\tilde{M}_{0}^{-1} & 0 \\
			0 & \tilde{N}_0^{-1} 
                        \end{bmatrix}
		\begin{bmatrix}
			I & 0 & -\tilde{u}_1\\
			0 & 1 & 0 \\
			0 & 0 & 1
		\end{bmatrix}.
\end{eqnarray}
It can be shown, using similar computations as in the proof
of Lemma 1, that
\begin{eqnarray*}
\tilde{N}_0^{-1} & =&  
\frac{1}{\tilde{\alpha}_1\tilde{\delta}_1 - \tilde{\beta}_1^2}
\begin{bmatrix}
\tilde{\delta}_1 &
-\tilde{\beta}_1 \\ -\tilde{\beta}_1
 & \tilde{\alpha}_1
\end{bmatrix} = 
\begin{bmatrix}
-\phi\lambda_1 & -s_1^Ty_1 - \phi \lambda_1 \\
-s_1^Ty_1 - \phi \lambda_1 & -s_1^Ty_1 - \phi \lambda_1 - y_1^TH_1y_1
\end{bmatrix}
.
\end{eqnarray*}
Substituting $\tilde{N}_0^{-1}$ into (\ref{eqn-biginv}) and simplifying yields
\begin{eqnarray}\label{eqn-big2}
\bar{M}_1^{-1}	=	\begin{bmatrix}
			\tilde{M}_{0}^{-1} & 0 & -\tilde{\Psi}_{0}^Ty_{1} \\
			0 & -\phi \lambda_1 & -s_1^Ty_1 - \phi \lambda_1  \\
			-y_1^T\tilde{\Psi}_{0}& -s_1^Ty_1 - \phi \lambda_1
& \tilde{u}_1^TM_0^{-1}\tilde{u}_1-s_1^Ty_1 - \phi \lambda_1 - y_1^TH_1y_1
		\end{bmatrix}.
\end{eqnarray}
Terms in the (3,3)-entry in the right side of (\ref{eqn-big2}) can be simplified using that
$\tilde{M}_{0}^{-1}\tilde{u}_1 = \tilde{\Psi}_{0}^Ty_1$  and
that $y_1^TH_1y_1 = y_1^TH_0y_1 + y_1^T\tilde{\Psi}_{0}\tilde{M}_{0}\tilde{\Psi}_{0}^Ty_1$ as follows:
\begin{eqnarray*}
	\tilde{u}_{1}^T\tilde{M}_{0}^{-1}\tilde{u}_1 - y_1^TH_1y_1
	&=& y_1^T\tilde{\Psi}_{0}\tilde{M}_{0}\tilde{\Psi}_{0}^Ty_1 -y_1^TH_1y_1\\
	&=& -y_1^TH_0y_1 + y_1^TH_1y_1 - y_1^TH_1y_1 \\
	&=& - y_1^TH_0y_1.
\end{eqnarray*}
Thus,
\begin{equation}\label{eqn-Mbarfinal}
\bar{M}_1^{-1}	=	\begin{bmatrix}
			\tilde{M}_{0}^{-1} & 0 & -\tilde{\Psi}_{0}^Ty_{1} \\
			0 & -\phi \lambda_1 & -s_1^Ty_1 - \phi \lambda_1  \\
			-y_1^T\tilde{\Psi}_{0}& -s_1^Ty_1 - \phi \lambda_1
& -s_1^Ty_1- \phi \lambda_1 - y_1^TH_0y_1
		\end{bmatrix}.
\end{equation}
%
Lemma 1, together with the substitution $\tilde{\Psi}_{0}^Ty_1 = \displaystyle
\begin{pmatrix}
	S_{0}^Ty_1 \\
	Y_{0}^TH_0y_1
\end{pmatrix}
$, 
gives that
\begin{eqnarray*}
&& \hspace{-.75cm}
\Pi_1^T \bar{M}_1^{-1} \Pi_1\\
	&=& 
	\Pi_j^T
	\begin{bmatrix}
			\tilde{M}_{0}^{-1} & 0 & -\tilde{\Psi}_{0}^Ty_{1} \\
			0 & -\phi \lambda_1 & -s_1^Ty_1 - \phi \lambda_1  \\
			-y_1^T\tilde{\Psi}_{0}& -s_1^Ty_1 - \phi \lambda_1
& -s_1^Ty_1- \phi \lambda_1 - y_1^TH_0y_1
		\end{bmatrix}
	\Pi_j \\
	&=& \tiny
	\left [
	\begin{tabular}{cc|ccccc}
		$-\phi\lambda_0$ & $0$ 
		& $-s_0^Ty_0-\phi\lambda_0$ & $-S_{0}^Ty_1$ \\
		$0$& $-\phi\lambda_1$ 
		& $0$ & $-s_1^Ty_1 - \phi \lambda_1$  \\\hline
		$-s_0^Ty_0-\phi\lambda_0$ & $0$ & $-s_0^Ty_0 - \phi\lambda_0 - y_0^TH_0y_0$ & $-Y_{0}^TH_0y_1$ \\
		$-y_1^TS_{0}$ & $-s_1^Ty_1 - \phi \lambda_1$ 
		& $-y_1^TH_0Y_{0}$ & $-s_1^Ty_1 - \phi \lambda_1 - y_1^TH_0y_1$
	\end{tabular}
	\right ] \\
	&=&	\tiny
	\left [
	\begin{tabular}{cc|ccccc}
		$-\phi \Lambda_{0}$ & $0$ 
		& $-R_{0} - D_{0} - \phi \Lambda_{0}$ & $-S_{0}^Ty_1$ \\
		$0$& $-\phi\lambda_1$ 
		& $0$ & $-s_1^Ty_1 - \phi \lambda_1$  \\\hline
		$-R_{0}^T - D_{0} - \phi \Lambda_{0} $ & $0$ & $-D_{0} - \phi\Lambda_1 - Y_{0}^TH_0Y_{0}$ & $-Y_{0}^TH_0y_1$ \\
		$-y_1^TS_{0}$ & $-s_1^Ty_1 - \phi \lambda_1$ 
		& $-y_1^TH_0Y_{0}$ & $-s_1^Ty_1 - \phi \lambda_1 - y_1^TH_0y_1$
	\end{tabular}
	\right ]\\[.3cm]
	&=& 
	\begin{bmatrix}
	-\phi \Lambda_{1} & -R_{1} - D_{1} - \phi \Lambda_{1} \\
	-R_{1}^T - D_{1} - \phi \Lambda_{1}  & -D_{1} - \phi\Lambda_1 - Y_{1}^TH_0Y_{1}	
	\end{bmatrix}\\[.35cm]
	&=&	
	\tilde{M}_{1}^{-1},
\end{eqnarray*}
proving the base case.

\medskip

The inductive step is similar to the base case.
For the inductive step, we assume the theorem holds for $k=1,\ldots j-1$
and now show it holds for $k=j$.  We begin, as before, setting
$k=j$ in \eqref{eqn-1paraminv} to obtain
\begin{equation}\label{eqn-induct}
	H_{j+1}= H_j + 
	 \begin{bmatrix}
                s_j & H_jy_j
        \end{bmatrix}
        \begin{bmatrix}
              \tilde{\alpha}_j & \tilde{\beta}_j \\ \tilde{\beta}_j & \tilde{\delta}_j
        \end{bmatrix}
        \begin{bmatrix}
                s_j^T \\
                y_j^TH_j
        \end{bmatrix},
\end{equation}
where
\begin{equation*}
	\tilde{\alpha}_j = \frac{1}{s_j^Ty_j}+ \Phi_j\frac{y_j^TH_jy_j}{(s_j^Ty_j)^2}, \quad
	\tilde{\beta}_j = -\frac{\Phi_j}{y_j^Ts_j}, \quad
\text{and}\quad	\tilde{\delta}_j = -\frac{1-\Phi_j}{y_j^TH_jy_j},
\end{equation*}
Using the induction hypothesis, (\ref{eqn-induct}) becomes
\begin{equation*} H_{j+1}=
	H_0 + \tilde{\Psi}_{j-1}\tilde{M}_{j-1}\tilde{\Psi}_{j-1} + 
	 \begin{bmatrix}
                s_j & H_jy_j
        \end{bmatrix}
        \begin{bmatrix}
              \tilde{\alpha}_j & \tilde{\beta}_j \\ \tilde{\beta}_j & \tilde{\delta}_j
        \end{bmatrix}
        \begin{bmatrix}
                s_j^T \\
                y_j^TH_j
        \end{bmatrix}.
\end{equation*}
Similar to the base case,
\begin{eqnarray*}
	&& \hspace{-1.75cm}
	\begin{bmatrix}
                s_j & H_jy_j
        \end{bmatrix}
        \begin{bmatrix}
              \tilde{\alpha}_j & \tilde{\beta}_j \\ \tilde{\beta}_j & \tilde{\delta}_j
        \end{bmatrix}
        \begin{bmatrix}
                s_j^T \\
                y_j^TH_j
        \end{bmatrix}
        \\
 	&=&
	\begin{bmatrix}
                s_j & H_0y_j +  \tilde{\Psi}_{j-1}\tilde{u}_{j}
        \end{bmatrix}
        \begin{bmatrix}
              \tilde{\alpha}_j & \tilde{\beta}_j \\ \tilde{\beta}_j & \tilde{\delta}_j
        \end{bmatrix}
        \begin{bmatrix}
                s_j^T \\
                y_j^TH_0 + \tilde{u}_{j}^T \tilde{\Psi}_{j-1}^T
        \end{bmatrix}\\
        &=&
        \begin{bmatrix}
        		\tilde{\Psi}_{j-1} & s_j & H_0y_j
        \end{bmatrix}
	\begin{bmatrix}
		\tilde{\delta}_j\tilde{u}_j\tilde{u}_j^T & \tilde{\beta}_j\tilde{u}_j 
			& \tilde{\delta}_j \tilde{u}_j \\
		\tilde{\beta}_j\tilde{u}_j^T & \tilde{\alpha}_j & \tilde{\beta}_j \\
		\tilde{\delta}_j\tilde{u}_j^T &\tilde{\beta}_j & \tilde{\delta}_j
	\end{bmatrix}
        \begin{bmatrix}
        		\tilde{\Psi}_{j-1}^T \\
		s_j^T \\
		y_j^TH_0
        \end{bmatrix},
\end{eqnarray*}
where $\tilde{u}_j=\tilde{M}_{j-1}\tilde{\Phi}^T_{j-1}y_j$.
With $\Pi_j$ defined as in \eqref{eqn-th1}, then
$
	\begin{bmatrix}
		\tilde{\Psi}_{j-1} & s_j & H_0y_j 
	\end{bmatrix}
	\Pi_j
	= \tilde{\Psi}_j,
$ and so
$
	H_j = H_0 +
	\tilde{\Psi}_j
	\Pi_j^T
\bar{M}_j	\Pi_j \tilde{\Psi}_j^T,
$
where
\begin{equation}\label{eqn-Mbarj}
\bar{M}_j\defined	\begin{bmatrix}
		\tilde{M}_{j-1}
			+\tilde{\delta}_j\tilde{u}_j\tilde{u}_j^T & \tilde{\beta}_j\tilde{u}_j 
			& \tilde{\delta}_j \tilde{u}_j \\
		\tilde{\beta}_j\tilde{u}_j^T & \tilde{\alpha}_j & \tilde{\beta}_j \\
		\tilde{\delta}_j\tilde{u}_j^T &\tilde{\beta}_j & \tilde{\delta}_j
	\end{bmatrix}.
\end{equation}

It remains to show that $\tilde{M}_j=\Pi_j^T\bar{M}_j\Pi_j$; however, for
simplicity, we instead show $\tilde{M}_j^{-1} = \Pi_j^T\bar{M}_j^{-1}\Pi_j$
using a similar argument as in the base case.
In particular, similar to \eqref{eqn-big2}, it can be shown that
$$
\bar{M}_j^{-1}=		\begin{bmatrix}
			\tilde{M}_{j-1}^{-1} & 0 & -\tilde{\Psi}_{j-1}^Ty_{j} \\
			0 & -\phi \lambda_j & -s_j^Ty_j - \phi \lambda_j  \\
			-y_j^T\tilde{\Psi}_{j-1}& -s_j^Ty_j - \phi \lambda_j & -s_j^Ty_j - \phi \lambda_j - y_j^TH_0y_j
		\end{bmatrix}.
$$
By the inductive hypothesis that (\ref{eqn-mtilde}) holds for $k=j-1$,
$$
	\tilde{M}_{j-1}^{-1} = 
	\begin{bmatrix}
	-\phi \Lambda_{j-1} & -R_{j-1} - D_{j-1} - \phi \Lambda_{j-1} \\
	-R_{j-1}^T - D_{j-1} - \phi \Lambda_{j-1}  & -D_{j-1} - \phi\Lambda_k - Y_{j-1}^TH_0Y_{j-1}	
	\end{bmatrix}
$$
and so,
\begin{eqnarray*}
	&& \hspace{-.5cm} \Pi_j^T \bar{M}_j^{-1}
	\Pi_j 	 \\
	&=& 
	\Pi_j^T
	\begin{bmatrix}
			\tilde{M}_{j-1}^{-1} & 0 & -\tilde{\Psi}_{j-1}^Ty_{j} \\
			0 & -\phi \lambda_j & -s_j^Ty_j - \phi \lambda_j  \\
			-y_j^T\tilde{\Psi}_{j-1}& -s_j^Ty_j - \phi \lambda_j & -s_j^Ty_j - \phi \lambda_j - y_j^TH_0y_j
	\end{bmatrix}	
	\Pi_j
	\\
	&=& \tiny
	\left [
	\begin{tabular}{cc|ccccc}
		$-\phi \Lambda_{j-1}$ & $0$ 
		& $-R_{j-1} - D_{j-1} - \phi \Lambda_{j-1}$ & $-S_{j-1}^Ty_j$ \\
		$0$& $-\phi\lambda_j$ 
		& $0$ & $-s_j^Ty_j \!-\! \phi \lambda_j$  \\\hline
		$-R_{j-1}^T \!\! - \! D_{j-1} \!\!-\! \phi \Lambda_{j-1} $ & $0$ & $-D_{j-1} \!\!-\! \phi\Lambda_{j-1}\! \!-\! Y_{j-1}^TH_0Y_{j-1}$ & $-Y_{j-1}^TH_0y_j$ \\
		$-y_j^TS_{j-1}$ & $-s_j^Ty_j \!-\! \phi \lambda_j$ 
		& $-y_j^TH_0Y_{j-1}$ & $-s_j^Ty_j \!-\! \phi \lambda_j \!-\! y_j^TH_0y_j$
	\end{tabular}
	\right ]\\[.3cm]
	&=& 
	\begin{bmatrix}
	-\phi \Lambda_{j} & -R_{j} - D_{j} - \phi \Lambda_{j} \\
	-R_{j}^T - D_{j} - \phi \Lambda_{j}  & -D_{j} - \phi\Lambda_j - Y_{j}^TH_0Y_{j}	
	\end{bmatrix}\\[.35cm]
	&=&	
	\tilde{M}_{j}^{-1},
\end{eqnarray*}
proving the induction step. $\square$

\bigskip

Using Theorem 1, Algorithm 1 computes $\tilde{M}_k$ and then uses the
compact formulation for the inverse to solve linear systems of the form
(\ref{eqn-basic}).  There are two parts to the \texttt{for} loop in
Algorithm 1: The first half of the loop is used to build products of the
form $s_j^TB_js_j$; the second half of the loop computes $\tilde{M}_k$,
making use of $s_j^TB_js_j$ to compute $\Phi_j$.  In the special case when
$\phi=0$ or $\phi=1$, $\Phi_j$ can be quickly computed as $\Phi_j=1$ and
$\Phi_j=0$, respectively.

\begin{algorithm}[H]
\caption{Computing $r=B_{k+1}^{-1}z$ using
    \eqref{eqn-H2}}
\tab Define $\phi$, $B_0$, and $H_0$;
\tab Form $M_0$ using (\ref{eqn-alt-form-convex}) and $\tilde{M}_0$ using \eqref{eqn-M0};
\tab Let $\Psi_0 = (B_0s_0 \,\, y_0 )$ and $\tilde{\Psi}_0 = (s_0 \,\, H_0y_0 )$;
\tab \FOR\ $j = 1:k$
\tabb  
\tabb $\%$ Part 1: Compute $s_j^TB_js_j$
\tabb $u_j \gets M_{j-1} (\Psi_{j-1}^Ts_j$);
\tabb $s_j^TB_js_j \gets s_j^TB_0s_j + ( s_j^T\Psi_{j-1})u_j$;
\tabb $\alpha_j \gets -(1-\phi)/ (s_j^TB_js_j$);
\tabb $\beta_j \gets -\phi/ (y_j^Ts_j$);
\tabb $\delta_j \gets  ( 1 + \phi (s_j^TB_js_j)/(y_j^Ts_j) ) / (y_j^Ts_j)$;
\tabb
$\displaystyle 
	M_j \ \gets \ 
	\Pi_j
	\begin{pmatrix}
			M_{j-1}+ \alpha_j u_ju_j^T & \alpha_j u_j & \beta_j u_j  \\
			\alpha_j u_j^T & \alpha_j & \beta_j \\
			\beta_j  u_j^T & \beta_j & \delta_j
	\end{pmatrix}
	\Pi_j^T
$, where $\Pi_j$ is as in (\ref{eqn-th1});
\tabb
\tabb $\%$Part 2: Compute $\tilde{M}_j$
\tabb $\tilde{u}_j \gets \tilde{M}_{j-1}(\tilde{\Psi}_{j-1}^Ty_j)$;
\tabb $y_j^TH_jy_j \gets y_j^TH_0y_j+ (y_j^T\tilde{\Psi}_{j-1})\tilde{u}_j$;
\tabb $\Phi_j = (1-\phi)(y_j^Ts_j)^2/((1-\phi)(y_j^Ts_j)^2 + \phi(y_j^TH_jy_j)(s_j^TB_js_j))$;
\tabb $\tilde{\alpha}_j \gets (1 + \Phi_j(y_j^TH_jy_j)/(y_j^Ts_j))/(y_j^Ts_j)$;
\tabb $\tilde{\beta}_j \gets -\Phi_j/(y_j^Ts_j)$;
\tabb $\tilde{\delta}_j \gets -(1-\Phi_j)/(y_j^TH_jy_j)$;
\tabb
$\displaystyle
	\tilde{M}_j \ \gets \ 
	\Pi_j^T
		\begin{pmatrix}
		\tilde{M}_{j-1}
			+\delta_j\tilde{u}_j\tilde{u}_j^T & \tilde{\beta}_j\tilde{u}_j
			& \tilde{\delta}_j \tilde{u}_j \\
		\tilde{\beta}_j\tilde{u}_j^T & \tilde{\alpha}_j & \tilde{\beta}_j  \\
		\tilde{\delta}_j\tilde{u}_j^T &\tilde{\beta}_j & \tilde{\delta}_j
		\end{pmatrix}
	\Pi_j
$, where $\Pi_j$ is as in \eqref{eqn-th1};
\tab \END\ \\
$r = H_0z+\tilde{\Psi}_k \tilde{M}_k \tilde{\Psi}_k^Tz$;
\end{algorithm}

In Algorithm 1, the computations for $u_j$ and
$\tilde{u}_j$ in lines 7 and 15, respectively, can be simplified by noting that
\begin{equation}\label{eqn-simple}
	\Psi_{j-1}^Ts_j =
	\begin{bmatrix}
		S_{j-1}^TB_0s_j \\
		Y_{j-1}^Ts_j
	\end{bmatrix}
\quad \text{and} \quad
	\tilde{\Psi}_{j-1}^Ty_j 
	=
	\begin{bmatrix}
		S_{j-1}^Ty_j \\
		Y_{j-1}^TH_0y_j
	\end{bmatrix}.
\end{equation}
In other words, $\Psi_{j-1}^Ts_j$ is a $2j$ vector whose first $j$ entries
are the first $j$ entries in the $(j+1)$th column of $S_k^TB_0S_k$ and
whose last $j$ entries are the first $j$ entries in the $(j+1)$th row of
$S_k^TY_k$; similarly, $\tilde{\Psi}_{j-1}^Ty_j$ is a vector whose first
$j$ entries are the first $j$ entries in the $(j+1)$th column of $S_k^TY_k$
and whose last $j$ entries are the first $j$th entries in the $(j+1)$th
column of $Y_k^TH_0Y_k$. Thus, computing $u_j$ and $\tilde{u}_j$ only
requires $2j(4j-1)$ flops additional flops after precomputing and storing
the quantities $S_k^TB_0S_k$, $Y_k^TH_0Y_k$, and $S_k^TY_k$.  Note that provided
$B_0$ is a scalar multiple of the identity matrix and $H_0 = B_0^{-1}$,
then forming $S_k^TB_0S_k$ and $Y_k^TH_0Y_k$ requires $(k+1)(k+2)n$ flops
each, and forming $S_k^TY_k$ requires $(k+1)^2(2n-1)$ flops.

Computing the scalar $s_j^TB_js_j$ in line 8 of Algorithm 1 requires performing one inner product between $\Psi_{j-1}^Ts_j$ and
$u_j$ and then summing this with $s_j^TB_0s_j$, resulting in an operation
cost of $4j$.  This is the same cost is associated with forming
$y_j^TH_jy_j$ in line 16.  The
scalar $\Phi_j$ requires ten flops to compute.  The
other scalars ($\alpha_j,
\beta_j, \delta_j, \tilde{\alpha}_j, \tilde{\beta}_j,$ and
$\tilde{\delta}_j$) require a total of 14 flops.
Forming the matrices $M_j$ and $\tilde{M}_j$ requires $24j^2 + 16j$ flops
each.  Thus, excluding the cost of
forming $S_k^TB_0S_k$, $S_k^TY_k$, and $Y_k^TH_0Y_k$,
the operation count for computing $\tilde{M}_k$ is given by
$$
	\sum_{j=1}^k 4j(4j-1)+8j+24+24j^2+16j = \frac{1}{3}\left ( 40k^3 + 90k^2+122k \right ).
$$
Finally, forming $r$ in line 23 requires
$k+1$ vector inner products and
$(2k+1)(n+k+1) + 2n$ additional flops.
Thus, the overall flop count for Algorithm 1 is
$(2k+1)(n+k+1) + 2n + (40k^3 + 90k^2+122k)/3 + (2n-1)(k+1),$
{which includes the cost of $k+1$ vector inner products.}

\bigskip

Table 1 presents the computational complexity of Algorithm 1 and the
algorithms in Appendices A-C for linear solves where the system matrix is
from the restricted Broyden class of updates.  The operational costs do not
include the computation for $S_k^TY_k, S_k^TB_0S_k$, and $Y_k^TH_0Y_k$,
which all can be efficiently updated as new quasi-Newton pairs are
computed.
\begin{table}[h]
\centering
\begin{tabular}{|c|c|}
\hline
Alg.\ &  Flop count (including vector inner products) \\ \hline
1 & $(2k+1)(n+k+1) + 2n + \frac{1}{3}(40k^3+90k^2+122k)+(2n-1)(k+1)$\\
3 & $4nk+3k+n + 2k(2n-1)$ \\
4 & $ \frac{1}{2}(k+1)(23kn+52n+12k+42)+(2n-1)(6(k+1)^2+7(k+1))$\\
6 & $(5n+3)(k+1)k+(10n+14)(k+1)+(2n-1)(3k+2)(k+1)$ \\
\hline
\end{tabular}
\caption{Computational complexity comparison of Algorithms 1, 3, 4, and 6
  when the system matrix is a member of the restricted Broyden class of quasi-Newton updates.}
\end{table}

\subsection{Efficient solves after updates}

To compute the factors in (\ref{eqn-H1}), it is necessary to
form the matrices $L_k$, $D_k$ and
$R_k$.  After a new quasi-Newton
pair has been computed, updating (\ref{eqn-H1}) can be done efficiently by storing $S_k^TS_k$, $Y_k^TY_k$
and $S_k^TY_k$.
In this case, to update $H_{k+1}$ we add a column to (and possibly delete a column
from) $S_k$ and $Y_k$; the corresponding changes can then be made in $S_k^TS_k$,
$Y_k^TY_k$ and $S_k^TY_k$ (see~\cite{ByrNS94} for more details).

In the practical implementation given in Algorithm 1, it is necessary to
compute $\Psi_{j-1}^Ts_j$ and $\tilde{\Psi}_{j-1}^Ty_j$ given by
\eqref{eqn-simple}, which can be formed using the stored quantities
$S_k^TS_k$, $Y_k^TY_k$, and $S_k^TY_k$.  To compute these quantities when a
new quasi-Newton pair of updates is obtained we apply the same strategy as
described above by adding (and possibly deleting) columns in $S_k$ and
$Y_k$, enabling $\Psi_{j-1}^Ts_j$ and $\tilde{\Psi}_{j-1}^Ty_j$ to be
updated efficiently. With these updates, the work required for
Algorithm 1 is significantly reduced.

\subsection{Compact representation of the inverse of an SR1 matrix}

In this subsection, we demonstrate that this same strategy can be applied to the \SR{} update.
The \emph{symmetric rank-one} (\SR) update is obtained by setting
$\phi = y_k^Ts_k/(y_k^Ts_k-s_k^TB_ks_k$; in this case,
\begin{eqnarray}\label{eqn-SR1}
	B_{k+1} &=& B_k  + \frac{1}
	{s_k^T(y_k  - B_k s_k ) }(y_k -B_k s_k )(y_k  - B_k s_k )^T.
\end{eqnarray}
The \SR{} update is a member of the Broyden class but not the restricted
Broyden class.  This update exhibits hereditary symmetric but not positive
definiteness.  This update is remarkable in that it is self-dual: 
Initializing $B_0^{-1}$ in place of $B_0$ and interchanging
$y_k$ and $s_k$ everywhere in (\ref{eqn-SR1}) yields a recursion
relation for $B_{k+1}^{-1}$.  Thus, like the \BFGS{} update, solving
linear systems with \SR{} matrices can be performed using vector
inner products (see  \ref{app:SR1} for details). 

\bigskip

The compact formulation
for the \SR{} update is given by $	B_{k+1}  = B_0  + \hat{\Psi}_k \hat{M}_k \hat{\Psi}_k^T,$
where 
\begin{equation}\label{eqn-alt-form-SR1}
	\hat{\Psi}_k = Y_k - B_0S_k \quad \text{and} \quad 
	\hat{M}_k  \ = \ (D_k + L_k + L_k^T - S_k^TB_0S_k)^{-1},
\end{equation}
where $\Psi_k \in \Re^{n \times (k+1)}$ and $M_k \in \Re^{(k+1)\times
  (k+1)}$ (see~\cite{ByrNS94}). 
The \SR{} update has the distinction of being the only rank-one update
in the Broyden class of updates.  It is for this reason that $\Psi_k$ 
 has only $k+1$ columns and $M_k$ is
a $(k+1)\times (k+1)$ matrix.  (In contrast,  the compact representation
for rank-two updates have twice as many columns and $M_k$ is twice
the size as in the \SR{} case.)
%
 
The compact formulation for the inverse of an \SR{}
matrix can be derived using the same strategy as for the 
restricted Broyden class.  The \SMW{} formula applied to
the compact formulation of an \SR{} with
$\Psi_k$ and $M_k$ defined as in (\ref{eqn-alt-form-SR1}) yields
\begin{eqnarray*}
	H_{k+1} 
	&=& 
	H_0 + H_0 (Y_k - B_0S_k)\left ( -M_k - \Psi_k^TH_0\Psi_k \right )^{-1}(Y_k-B_0S_k)^TH_0.
\end{eqnarray*}
Substituting in for $M_{k}$ using (\ref{eqn-alt-form-SR1}) together with the
identity 
$$
	\Psi_k^TH_0\Psi_k = Y_k^TH_0Y_k - S_k^TY_k - Y_k^TS_k + S_k^TB_0S_k,
$$
yields
\begin{eqnarray*}
H_{k+1}
	&=&
	H_0 +  (H_0 Y_k - S_k)
	 \big ( 
	 - (D_k+L_k+L_k^T-S_k^TB_0S_k) \\
	 && \qquad \qquad 
	 - (Y_k^TH_0Y_k - S_k^TY_k - Y_k^TS_k + S_k^TB_0S_k)
	 \big )^{-1}
	 (H_0Y_k - S_k)^T\\
	 &=& 
	 H_0 + (S_k- H_0Y_k)^T( D_k + R_k + R_k^T -Y_k^TH_0Y_k )^{-1}(S_k - H_0Y_k)^T.
\end{eqnarray*}
In other words, 
the compact formulation for the inverse of an \SR{} matrix
is given by
\begin{equation}\label{eqn-H3}
	H_{k+1}= 
	H_0 + \tilde{\Psi}_k \tilde{M}_k \tilde{\Psi}_k^T \quad 
\text{where}\quad \tilde{\Psi}_k \defined [S_k-H_0Y_k],
\end{equation}
and 
$
\tilde{M}_k \defined 
(D_k + R_k + R_k^T-Y_k^TH_0Y_k)^{-1}.
$  Algorithm 2 details how to solve a linear system defined by an \LSR{} matrix
via the compact representation for its inverse.

\begin{algorithm}
\caption{Computing $r = B_{k+1}^{-1}z$ using \eqref{eqn-H3} when $B_{k+1}$ is an SR1 matrix}
\tab Define $H_0, S_k$ and $Y_k$;
\tab Form $S_k^TY_k = L_k + D_k + R_k$;
\tab Form $\tilde{\Psi}_{k} = S_k - H_0Y_k$;
\tab $\tilde{M}_k \leftarrow (D_k + R_k + R_k^T - Y_k^TH_0Y_k)^{-1}$;
\tab $r = H_0z + \tilde{\Psi}_k\tilde{M}_k\tilde{\Psi}_k^Tz$;
\end{algorithm}

Assuming $S_k^TY_k$ is precomputed, Algorithm 2 requires 
$(2k+1)(n+k+1) + 2n + (2n-1)(k+1)$ flops, which includes
$(k+1)$ vector inner products to compute $\tilde{\Psi}_k^Tz$ in line 5.
Note that this count does not include the cost of inverting the
  $(k+1) \times (k+1)$ matrix, $\tilde{M}_k$, in line 4 since when $k \ll
  n$, this cost is significantly smaller than the dominant costs 
 for Algorithm 2.

\bigskip 

Table 2 presents the computational complexity of Algorithm 2
and Algorithm 8 in \ref{app:SR1} for solving SR1 linear systems.
The operational costs do not include the computation
for $S_k^TY_k, S_k^TB_0S_k$, and $Y_k^TH_0Y_k$, which all can be
efficiently updated as new quasi-Newton pairs are computed.

\begin{table}[h]
\centering
\begin{tabular}{|c|c|}
\hline
Algorithm & Flop count (including vector inner products)\\ \hline
2 & $(2k+1)(n+k+1) + 2n + (2n-1)(k+1)$\\
8   & $\frac12(k+1)(k(2n+1)+2(3n+1))+2n + \frac12(2n-1)(k+2)(k+1)$ \\
\hline
\end{tabular}

\caption{Computational complexity comparison of Algorithms 2 and 8 when
  the system matrix is an SR1 matrix.  Note that the computational cost for Algorithm
  2 does not include the cost of inverting a $(k+1) \times (k+1)$ matrix, where $k \ll n$.}
\end{table}

\subsection{Computing condition numbers}
Provided the initial approximate Hessian $B_0$ is taken to be a scalar
multiple of the identity, (i.e., $B_0=\gamma I$), it is possible to compute
the condition number of $B_{k+1}$ in (\ref{eqn-basic}).  To do this, we
consider the condition number of $H_{k+1}=B_{k+1}^{-1}$.  The eigenvalues
of $H_{k+1}$ can be obtained from the compact representation of $H_{k+1}$
together with the techniques from~\cite{ErwM15,Burdakov13}.  For
completeness, we review this approach below.

Suppose $B_{k+1}$ is a member of the restricted Broyden class.  The
eigenvalues of $H_{k+1}$ can be computed using the compact formulation for
$H_{k+1}$ given in (\ref{eqn-H2}):
$$
	H_{k+1}=H_0 + \tilde{\Psi}_k \tilde{M}_k \tilde{\Psi}_k^T.
$$
Letting $\tilde{\Psi}_k=QR$ be the \QR{} decomposition of $\tilde{\Psi}_k$
where $Q\in\Re^{n\times n}$ is an orthogonal matrix and $R\in\Re^{n\times
  2(k+1)}$ is an upper triangular matrix.  Then, 
\begin{eqnarray}
	H_{k+1} &= &
	H_0 + \tilde{\Psi}_k \tilde{M}_k \tilde{\Psi}_k^T\\  \nonumber
	& = & H_0 + QR\tilde{M}_k R^TQ^T \label{eqn-eig1}.
\end{eqnarray}
The matrix $R\tilde{M}_kR^T$ is a real symmetric $n\times n$ matrix.
However, since $\tilde{\Psi}_k\in\Re^{n\times 2(k+1)}$ has at most
rank $2(k+1)$, then $R$ can be written in the form
$$
	R=\begin{pmatrix} R_1 \\  0 \end{pmatrix},
$$
where $R_1\in 2(k+1) \times 2(k+1)$.  Thus,
$$
	R\tilde{M}_kR^T=\begin{pmatrix}
	R_1 \\ 0 \end{pmatrix} \tilde{M}_k \begin{pmatrix} R_1^T & 0\end{pmatrix}=
	\begin{pmatrix}
		R_1\tilde{M}_kR_1^T & 0 \\ 0 & 0
	\end{pmatrix}.
$$
Since $R_1\tilde{M}_kR_1^T\in\Re^{2(k+1)\times 2(k+1)}$, its spectral
decomposition can be computed explicitly.  Letting $V_1D_1V_1^T$
is the spectral decomposition of $R_1\tilde{M}_kR_1^T$ and
substituting into (\ref{eqn-eig1}) yields:
$$
	H_{k+1} =  H_0+QVDV^TQ^T
$$
where 
$$
	V\defined 
	\begin{pmatrix}
		V_1 & 0 \\ 0 & 0 
	\end{pmatrix}\in\Re^{n\times n} 
	\quad
	\text{and} \quad 
	D\defined 
	\begin{pmatrix} 
		D_1 & 0 \\ 0 & 0 
	\end{pmatrix} \in\Re^{n\times n}.
$$
Thus,
$$
	H_{k+1} =  H_0+QVDV^TQ^T = QV(\gamma^{-1} I + D )V^TQ^T,
$$
giving the spectral decomposition of $H_{k+1}$.  
In particular, the matrix $H_{k+1}$ has an eigenvalue of $\gamma^{-1}$
with multiplicity $n-2(k+1)$ and $2(k+1)$ eigenvalues given by
$\gamma^{-1}+d_i$ for $1\le i \le 2(k+1)$ where $d_i$ denotes
the $i$th diagonal element of $D_1$.  Thus, the condition number
of $H_{k+1}$ can be computed as the ratio of the largest eigenvalue
to the smallest eigenvalue (in magnitude).  The condition number
of $B_{k+1}$ is the reciprocal of this value.  

When $B_{k+1}$ a quasi-Newton matrix generated by \SR{} updates,
the procedure is similar except $\tilde{\Psi}_k$ has half as many columns
resulting in $(k+1)$ eigenvalues given by
$\gamma^{-1}+d_i$ for $1\le i \le (k+1)$ and $n-(k+1)$ eigenvalues
of $\gamma^{-1}$.
After a new quasi-Newton pair is computed it is possible to update
the \QR{} factorization (for details, see \cite{ErwM15}).

\section{Numerical experiments}

In this section, we demonstrate the accuracy of the proposed method for solving
quasi-Newton equations.  We solve the following linear system:
\begin{equation}\label{eqn-numex}
	B_{k+1}p_{k+1}= -g_{k+1},
\end{equation}
where $B_{k+1}$ is a member of the restricted Broyden
class of updates or an SR1 matrix. 
For members of the restricted Broyden class of updates, we considered three
values of $\phi$: $0.00, 0.50$, and $0.99$ and compare the following
methods:
\begin{enumerate}
	\item``Two-loop recursion'' for the BFGS case ($\phi = 0.0$) (Algorithm 3)
	\item Recursive SMW approach (Algorithm 4)
	\item Recursion to compute products with $B_{k+1}^{-1}  = H_{k+1}$ (Algorithm 6)
	\item Compact inverses formulation  (Algorithm 1)
\end{enumerate}
For SR1 matrices, we compare the
following methods:
\begin{enumerate}
	\item The self-duality method (Algorithm 8)
	\item Compact inverses formulation (Algorithm 2)
\end{enumerate}

We consider problem sizes ($n$) ranging from $10,000$ to $1,000,000$.
We present the norms of the relative residuals:
$$
	\frac{\| B_{k+1} p_{k+1} - (- g_{k+1})\|_2}{\| g_{k+1} \|_2},
$$
and the time (in seconds) needed to compute each solution.  
The number of limited-memory updates was set to five.

We simulate the first five iterations of an unconstrained line-search method.  
We generate the initial point $x_0$ and the gradients $g_j = \nabla f(x_j)$, for $0 \le j \le 5$, randomly and 
compute the subsequent iterations
$$
	x_{j+1} = x_j - \alpha_j H_{j}g_j, \qquad \text{for $0 \le j \le 5$},
$$
where $H_j=B_j^{-1}$ is the inverse of the quasi-Newton matrix.

\bigskip
\noindent
{\bf Results.} 
We ran each algorithm ten times for each value of $\phi$ $(0.00, 0.50$, and
$0.99$) and for each $n$ ($10,000$, $50,000$, $100,000$ and $1,000,000$).
The computational time results are shown in semi-log plots in Figure 1.
Representative relative error results are presented in Tables 3--6.  All
algorithms were able to solve the linear systems to high accuracy.  
For the restricted Broyden class of matrices, the proposed compact inverse
algorithm (Algorithm 1)
outperformed all methods for all sizes of systems except
for the \LBFGS{} case.  In this case, Algorithm 1 is 
comparable to the
the two-loop recursion (Algorithm 3), which is specific to the 
\LBFGS{} update.  
(When
$n=1,000,000$, Algorithm 1 was slightly
more efficient than Algorithm 3.)
For
the \SR{} update,
the compact inverse formulation outperformed the algorithm based on
self-duality (Algorithm 8).

\begin{figure}[p]
\label{fig:results}
\begin{tabular}{cc}
 \hspace{-.25cm} \includegraphics[width=6.3cm]{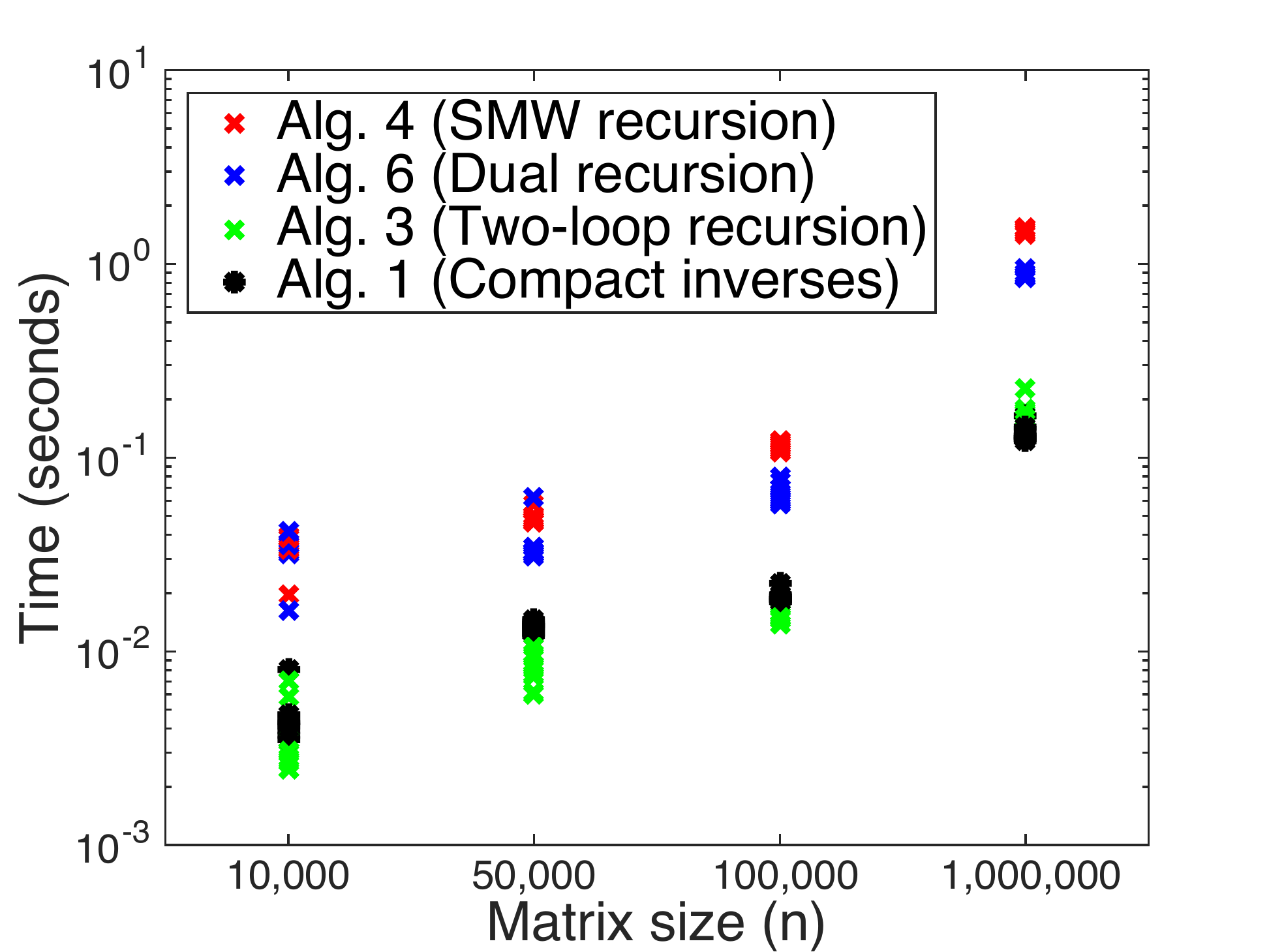} &  
 \hspace{-.75cm} \includegraphics[width=6.3cm]{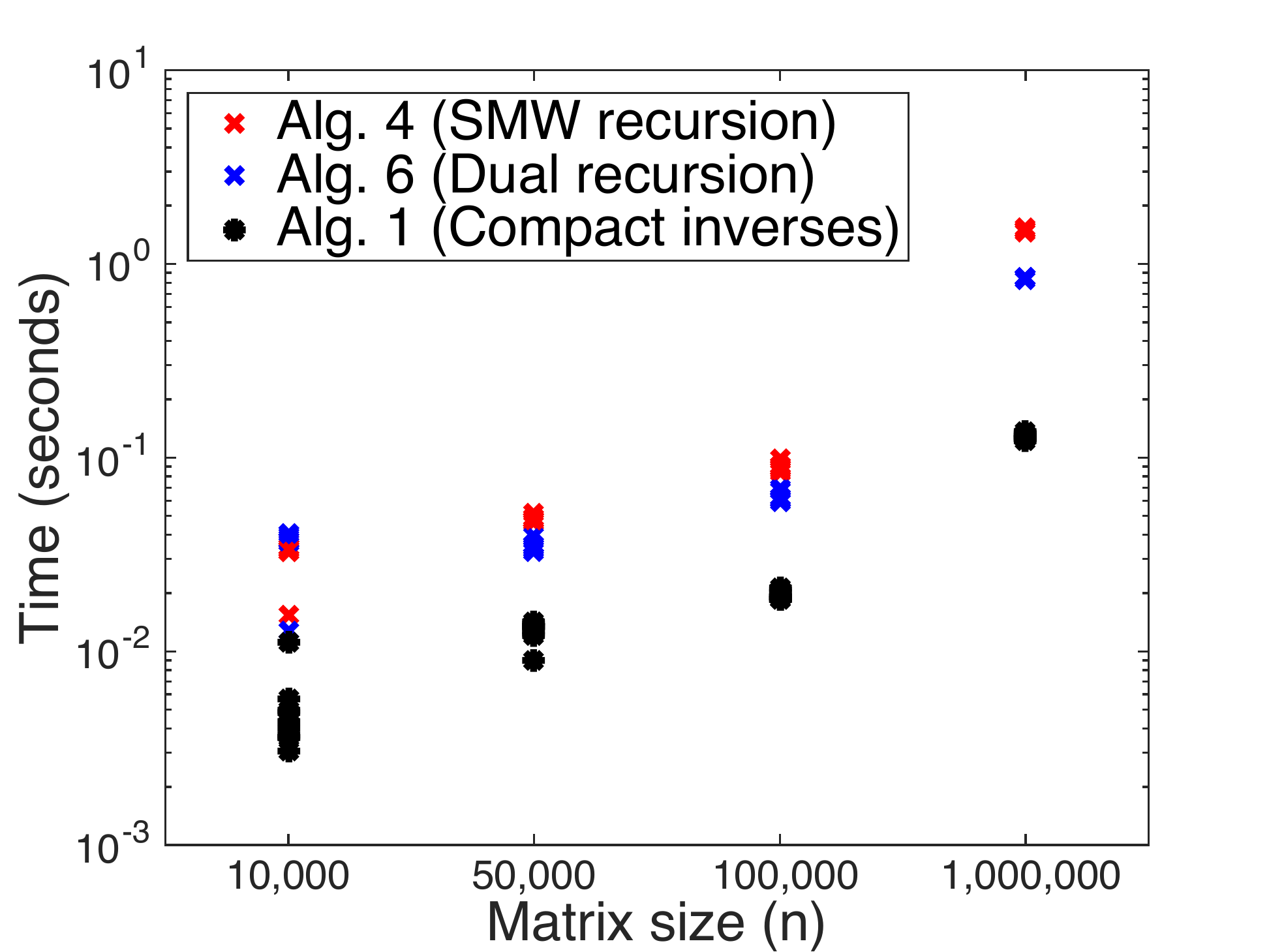} \\
(a) $\phi = 0$ (BFGS) & (b) $\phi = 0.5$ \ \ \  \\[.05cm]
 \hspace{-.25cm} \includegraphics[width=6.3cm]{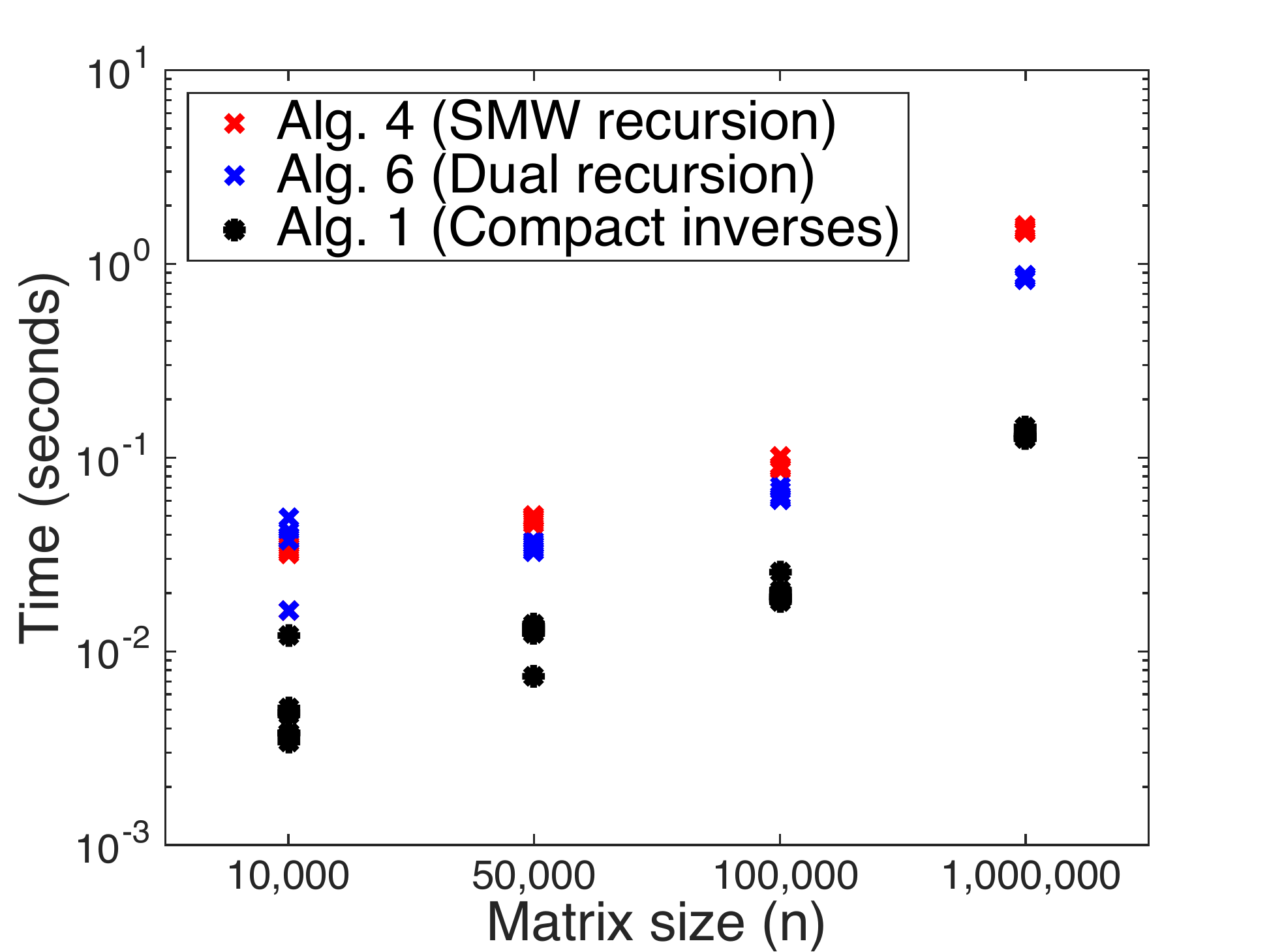} & 
 \hspace{-.75cm} \includegraphics[width=6.3cm]{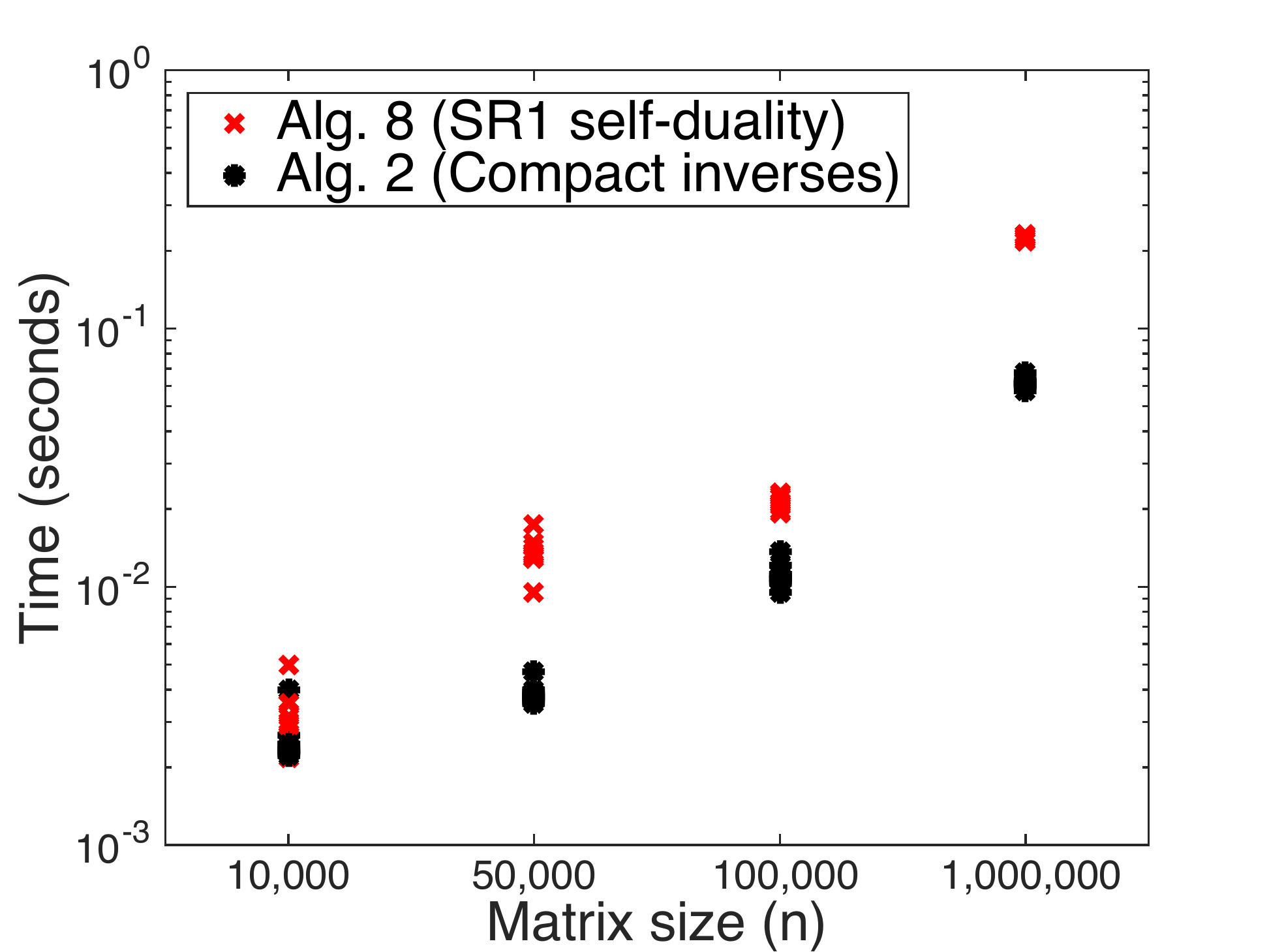} \\
(c) $\phi = 0.99$ & (d) SR1 \ \ \ 
\end{tabular}
\caption{Semi-log plots of the computational times (in seconds) for ten runs
  of each of the algorithms discussed in this paper. 
In (a),
  (b), and (c), the system matrix is 
a member of the restricted Broyden class of updates with
  $\phi = 0, 0.5,$ and $0.99$, respectively.  In (d), 
the system matrix is an \SR{} matrix.
The proposed method
  using compact inverses generally outperforms the other methods. 
  Note that when
  $\phi = 0$, our proposed method is competitive with the ``two-loop
  recursion'', which is specific to the BFGS update.}
\end{figure}

\setlength\tabcolsep{1.5mm}
\begin{table}[th]

\caption{Relative error for Broyden class matrices with $\phi = 0.00$ (\BFGS{}). 
} 
\begin{center}
\begin{tabular}{|r|>{\centering}m{2.2cm}|>{\centering}m{2.2cm}|>{\centering}m{2.2cm}|m{2.2cm}|}
\hline
& Two-Loop & Recursive & Recursive & \centering {Compact} \tabularnewline[-.2cm]
$n$ \quad \ \ \ & Recursion & SMW & $H_{k+1}$ & \centering {Inverses} \tabularnewline[-.2cm]
& (Alg.\ 3) & (Alg.\ 4) & (Alg.\ 6) & \centering {(Alg.\ 1)} \tabularnewline
\hline
\hline \texttt{ 10,000} & \texttt{4.88e-16} & \texttt{4.63e-16} & \texttt{5.37e-16} & \centering \texttt{3.59e-16} \tabularnewline
\hline \texttt{ 50,000} & \texttt{2.93e-16} & \texttt{4.44e-16} & \texttt{3.61e-16} & \centering \texttt{4.20e-16} \tabularnewline
\hline \texttt{100,000} & \texttt{6.64e-16} & \texttt{3.74e-16} & \texttt{6.16e-16} & \centering \texttt{3.81e-16} \tabularnewline
\hline \texttt{1,000,000} & \texttt{1.47e-15} & \texttt{1.46e-15} & \texttt{1.45e-15} & \centering \texttt{1.51e-15} \tabularnewline
\hline
\end{tabular} 
\end{center}
\end{table}

%
%

\setlength\tabcolsep{1.5mm}
\begin{table}[th]

\caption{Relative error for Broyden class matrices with $\phi = 0.50$.
} 
\begin{center}
\begin{tabular}{|r|>{\centering}m{2.8cm}|>{\centering}m{2.8cm}|m{2.8cm}|}
\hline
\multirow{2}{*}{$n$ \quad \ \ \ } & Recursive SMW & Recursive $H_{k+1}$ & \centering {Compact Inverses} \tabularnewline[-.2cm]
& (Alg.\ 4) & (Alg.\ 6) & \centering {(Alg.\ 1)} \tabularnewline
\hline
\hline \texttt{ 10,000} \ & \texttt{9.90e-16} & \texttt{8.37e-16} & \centering \texttt{8.15e-16} \tabularnewline
\hline \texttt{ 50,000} \ & \texttt{4.25e-16} & \texttt{6.80e-16} & \centering \texttt{5.82e-15} \tabularnewline
\hline \texttt{100,000} \ & \texttt{7.00e-16} & \texttt{6.31e-16} & \centering \texttt{9.14e-16} \tabularnewline
\hline \ \ \texttt{1,000,000} \ & \texttt{2.77e-16} & \texttt{2.40e-16} & \centering \texttt{3.56e-16} \tabularnewline
\hline
\end{tabular} 
\end{center}
\end{table}

%

\setlength\tabcolsep{1.5mm}
\begin{table}[th]

\caption{Relative error for Broyden class matrices with $\phi = 0.99$.
} 
\begin{center}
\begin{tabular}{|r|>{\centering}m{2.8cm}|>{\centering}m{2.8cm}|m{2.8cm}|}
\hline
\multirow{2}{*}{$n$ \quad \ \ \ } & Recursive SMW & Recursive $H_{k+1}$ & \centering {Compact Inverses} \tabularnewline[-.2cm]
& (Alg.\ 4) & (Alg.\ 6) & \centering {(Alg.\ 1)} \tabularnewline
\hline
\hline \texttt{ 10,000} \ & \texttt{8.33e-16} & \texttt{1.59e-15} & \centering \texttt{1.63e-15} \tabularnewline
\hline \texttt{ 50,000} \ & \texttt{8.45e-15} & \texttt{1.21e-14} & \centering \texttt{3.88e-15} \tabularnewline
\hline \texttt{100,000} \ & \texttt{7.28e-14} & \texttt{6.67e-14} & \centering \texttt{2.67e-14} \tabularnewline
\hline \ \ \texttt{1,000,000} \ & \texttt{2.59e-15} & \texttt{1.80e-15} & \centering \texttt{3.29e-15} \tabularnewline
\hline
\end{tabular} 
\end{center}
\end{table}

\setlength\tabcolsep{1.5mm}
\begin{table}[th]

\caption{Relative error for SR1 matrices. 
} 
\begin{center}
\begin{tabular}{|r|>{\centering}m{3.2cm}|>{\centering}m{3.2cm}|>{\centering}m{2.2cm}|m{2.2cm}|}
\hline
\multirow{2}{*}{$n$ \quad \ \ }  & Self-Duality & Compact Inverses \tabularnewline[-.2cm]
& (Alg.\ 8) & (Alg.\ 2)\tabularnewline
\hline
\hline \ \  \ \ \texttt{ 10,000} 
& \texttt{1.98e-15} 
& \centering \texttt{6.10e-15} 
\tabularnewline
\hline \ \ \  \ \texttt{ 50,000} 
& \texttt{2.24e-14} 
& \centering \texttt{7.57e-14} 
\tabularnewline
\hline \ \ \ \ \texttt{100,000}
& \texttt{5.07e-14} 
& \centering \texttt{6.44e-14} 
\tabularnewline
\hline \texttt{1,000,000} 
& \texttt{8.67e-13} 
& \centering \texttt{2.26e-12} 
\tabularnewline
\hline
\end{tabular} 
\end{center}
\end{table}

\section{Concluding Remarks}

We derived the compact formulation for members of the restricted Broyden
class and the \SR{} update.  With this compact formulation, we showed how
to solve linear systems defined by limited-memory quasi-Newton matrices.
Numerical results suggest that this proposed approach is efficient and
accurate.  This approach has two distinct advantages over existing
procedures for solving limited-memory quasi-Newton systems.  First, there
is a natural way to use the compact formulation for the inverse to obtain
the condition number of the linear system.  Second, when a new quasi-Newton
pair is computed, computational savings can be achieved by simple updates
to the matrix factors in the compact formulation.

Future work includes integrating this linear solver inside line-search
and trust-region methods for large-scale optimization.

\section{Acknowledgements}

The authors would like to thank Tammy Kolda for initial conversations on this subject.

\appendix

\section{Two-loop recursion}\label{app:two-loop}
The Broyden-Fletcher-Goldfarb-Shanno (\BFGS) update  is obtained by
setting $\phi=0$ in  \eqref{eqn-broyden}.  In this case, $B_{k+1}$ simplifies to
\begin{equation}\label{eqn-BFGS}
	B_{k+1} = B_k - \frac{1}{s_k^TB_ks_k}B_ks_ks_k^TB_k+ \frac{1}{y_k^Ts_k}y_ky_k^T.
\end{equation}
The inverse of the \BFGS{} matrix $B_{k+1}$ is given by
$$
	B_{k+1}^{-1} = 
	\left ( 
	I - \frac{s_ky_k^T}{y_k^Ts_k} 
	\right )
	B_k^{-1}
	\left (
	I - \frac{y_ks_k^T}{y_k^Ts_k}
	\right )
	+
	\frac{s_ks_k^T}{y_k^Ts_k},
$$
which can be written recursively as
\begin{eqnarray}
  	B^{-1}_{k+1} &=&  (V_{k}^T\cdots V_{0}^T ) B_0^{-1} (V_{0}\cdots V_{k}) +
  	\frac{1}{y_{0}^Ts_{0}} (V_{k}^T\cdots V_{1}^T)s_0s_0^T(V_1\cdots V_{k}) \nonumber\\
  	&& +
  	\frac{1}{y_{1}^Ts_{1}} (V_{k}^T\cdots V_{2}^T)s_1s_1^T(V_2\cdots V_{k})
  	+ \cdots + \frac{1}{y_{k}^Ts_{k}} s_{k}s_{k}^T, \label{eqn-lbfgs2}
\end{eqnarray}
where $V_i=I-\frac{1}{y_i^Ts_i}y_is_i^T$ (see, e.g.,~\cite[p.177--178]{NocW06}).  
Solving \eqref{eqn-basic} can be done efficiently using the
well-known two-loop recursion~\cite{Noc80}.

\bigskip

\begin{algorithm}

\caption{Two-loop recursion to compute $r=B_k^{-1}z$ when $B_{k}$ is a BFGS matrix}\label{alg-recursion}
  \tab $q\leftarrow z$;
  \tab \FOR\ $i=k-1,\dots,0$
  \tabb $\alpha_i\leftarrow(s_i^Tq)/(y_i^Ts_i)$;
  \ $q \leftarrow q-\alpha_iy_i$;
  \tab \END\ 
  \tab $r\leftarrow B_0^{-1} q$;
  \tab \FOR\ $i=0,\ldots, k-1$
  \tabb  $\beta \leftarrow (y_i^Tr)/(y_i^Ts_i)$;
  \  $r \leftarrow r+(\alpha_i-\beta)s_i$:
  \tab \END\
\end{algorithm}
Assuming that the inner products $y_i^Ts_i$ can be precomputed and $B_0$ is
a scalar multiple of the identity matrix, then the total operation count
for Algorithm 3 is $4nk + 3k + n + 2k(2n-1)$ flops, which includes $2k$
vector inner products.



\section{SMW recursion}\label{app:SMW-recursion}

In the symmetric case,
the \SMW{} formula for a rank-one change is given by
\begin{equation}\label{eqn-smw}
(A+\alpha uu^T)^{-1} = A^{-1} - \frac{\alpha}{1+\alpha u^TA^{-1}u} A^{-1}uu^TA^{-1}.
\end{equation}
We now show how (\ref{eqn-smw}) can be used to compute the inverse of 
$B_{k+1}$.  First, for $0\le j \le k$, define 
\begin{equation}\label{eqn-alphau}
\begin{array}{rlllrll}
\alpha_{3j} &=& (y_j^Ts_j )^{-1},& \hspace*{1cm}  u_{3j}&=& y_j, \\
\alpha_{3j+1} &=&  \phi (s_j^T B_j   s_j ), & \hspace*{1cm} u_{3j+1} &=&
	\displaystyle \frac{y_j}{y_j^Ts_j} - \frac{B_j   s_j}{s_j^TB_j   s_j},  \\
\alpha_{3j+2} &=& -(s_j^TB_j s_j )^{-1}, &\hspace*{1cm} u_{3j+2}&=& B_j  s_j, 
\end{array}
\end{equation}
Letting 
$$C_0=B_0 \quad \text{and} \quad U_i=\alpha_iu_iu_i^T
$$ 
for $0\le i \le 3k$, we define the matrices
\begin{equation}\label{eqn-Ck}
	C_1 = C_0+U_0, \quad C_2=C_1+U_1,  
	\quad \ldots, \quad 
	C_{3(k+1)}=C_{3(k+1)-1}+U_{3(k+1)-1} .
\end{equation}
By construction,
$C_{3(k+1)}=B_{k+1}$.  It can be shown that each $C_i$ is positive definite
and so $u_i^TC_i^{-1}u_i>0$ for each $i$.  For the restricted Broyden
class, $\alpha_{3j}$ and $\alpha_{3j+1}$, $j=0,\ldots, k$  are both
strictly positive.  Thus, the denominator in (\ref{eqn-smw}) is strictly
positive for the rank-one changes associated the $3j$ and $3j+1$,
$j=0,\ldots k$, and thus, the \SMW{} formula is well-defined for these
rank-one updates.  Since $C_{3j}=B_j$ and $B_j$ is positive definite, the
rank-one update associated with $j+2$, $j=0,\ldots k$, must result in a
positive-definite matrix; in other words, (\ref{eqn-smw}) is well
defined. Thus, solving the system $B_{k+1}r = z$ is equivalent to computing
$r = C_{3(k+1)}^{-1}z$.  Apply the \SMW{} formula in (\ref{eqn-smw}) to
obtain the inverse of $C_{i+1}$ from $C_i^{-1}$, we obtain
\begin{equation}\label{eqn-cinv}
	C_{i+1}^{-1}z = C_i^{-1}z-\frac{\alpha_i}{1+\alpha_i u_{i}^TC_i^{-1}u_{i}}C_i^{-1}u_iu_i^TC_i^{-1}z,
\end{equation}
for $0\le i < 3(k+1).$
Recursively applying (\ref{eqn-cinv}) to $C_i^{-1}z$, we obtain that 
\begin{equation}\label{eqn-cinv1b}
	C_{i+1}^{-1}z = C_0^{-1}z-\sum_{j=0}^{i}\frac{\alpha_j}{1+\alpha_j u_{j}^TC_j^{-1}u_{j}}C_j^{-1}u_ju_j^TC_j^{-1}z,
\end{equation}
and more importantly,
\begin{equation}\label{eqn-cinv2}
	C_{3(k+1)}^{-1}z = C_0^{-1}z-\sum_{j=0}^{3(k+1)-1}\frac{\alpha_j}{1+\alpha_j u_{j}^TC_j^{-1}u_{j}}C_j^{-1}u_ju_j^TC_j^{-1}z.
\end{equation}
Computing (\ref{eqn-cinv2}) can be simplified by defining the following
two quantities:
\begin{equation}\label{eqn-tauP}
	\tau_j = \frac{\alpha_j}{1+\alpha_j u_{j}^TC_j^{-1}u_{j}} \quad 
	\text{and} \quad 
	p_j = C_j^{-1}u_j.
\end{equation}
Substituting (\ref{eqn-tauP}) into (\ref{eqn-cinv2}) yields
\begin{equation}\label{eqn-cinv3}
	C_{3(k+1)}^{-1}z = C_0^{-1}z-\sum_{j=0}^{3(k+1)-1}\tau_j(p_j^Tz)p_j,
\end{equation}
where $\tau_j = \alpha_j(1+\alpha_j p_{j}^Tu_{j})^{-1}.$
The advantage of representation is that computing $C_{3(k+1)}^{-1}z$ requires
only vector inner products.  Moreover, the vectors $p_j$ can be computed
by evaluating (\ref{eqn-cinv1b}) at $z=u_j$ for $i=j-1$; that is,
$$p_j=C_j^{-1}u_j = C_0^{-1} u_j - \sum_{i=0}^{j-1}\tau_i(p_i^Tu_j)p_i.$$

In Algorithm 4, we present the recursion to solve the linear system
$B_{k+1} r =z$ using the \SMW{} formula (see related methods in~{\cite{ErwM12,Mill81}).  
We assume $B_0$ is an easily-invertible initial matrix. 

\begin{algorithm}[H]
\caption{Computing $r=B_{k+1}^{-1} z$}\label{alg-Hrec}
  \tab $r\leftarrow C_0^{-1}z$;
  \tab \FOR\ $j=0,\ldots,3(k+1)-1$
  \tabb \IF\ $\mod(j,3)=0$ 
  \tabbb $i\leftarrow j/3$; 
  \  $u\leftarrow y_{i}$; 
  \  $\alpha\gets y_{i}^Ts_{i}$;
  \tabb \ELSEIF\ $\mod(j,3)=1$
  \tabbb $i\leftarrow (j-1)/3$; \ $u \gets (y_{i})/(y_{i}^Ts_{i}) - (B_{i}s_{i})/(s_{i}^TB_{i}s_{i})$;
  \  $\alpha\gets \phi(s_{i}^TB_{i}s_{i})$;
  \tabb  \ELSE\ 
  \tabbb $i\leftarrow (j-2)/3$;
  \  $u\gets B_{i}s_{i}$;
  \  $\alpha\gets -(s_{i}^TB_{i}s_{i})^{-1}$;
  \tabb \END\ 
  \tabb $p_j\gets C_0^{-1}u$;
  \tabb \FOR\ $l=0,\ldots, j$
  \tabbb $p_j \gets p_j-\tau_l (p_l^Tu)p_l$;
  \tabb \END\
  \tabb $\tau_j\gets \alpha/(1+\alpha p_j^Tu)$;
  \tabb $r\leftarrow r - \tau_j(p_j^Tz)p_j$;
  \tab \END\
\end{algorithm}

\bigskip

At each iteration, Algorithm 4 computes $\alpha$ and $u$, which
are defined in \eqref{eqn-alphau} and require matrix-vector products
involving the matrices $B_i$ for $0 \le i \le k$.  The main difficulty in
computing $\alpha$ and $u$ is the computation of matrix-vector products
with the matrices $B_i$ for each $i$.  Note that if we are able to form
$u_{3j+2}=B_js_j$, then we are able to compute all
other terms that use $B_js_j$ in \eqref{eqn-alphau}.  In what follows, we
show how to compute $u_{3j+2}$ without storing $B_i$ for $0\le i \le k$.
This idea is based on~\cite[Procedure 7.6]{NocW06}.  

We begin by writing $C_{3(k+1)}$ in \eqref{eqn-Ck} as
follows:
\begin{eqnarray*}
	C_{3(k+1)} 
		&=& B_0 	+ \sum_{j=0}^k \alpha_{3j}u_{3j}u_{3j}^T 
				+ \alpha_{3j+1}u_{3j+1}u_{3j+1}^T 
				+ \alpha_{3j+2}u_{3j+2}u_{3j+2}^T,
\end{eqnarray*}
where $\alpha_i$ and $u_i$ are defined by \eqref{eqn-alphau}.  
Since $B_j = C_{3j}$, then
\begin{eqnarray*}
	u_{3j+2} 	&=& B_js_j \\
			&=&  \left ( B_0
				+ \sum_{i=0}^{j-1} \alpha_{3i}u_{3i}u_{3i}^T 
				+ \alpha_{3i+1}u_{3i+1}u_{3i+1}^T 
				+ \alpha_{3i+2}u_{3i+2}u_{3i+2}^T \right )s_j \\
			&=& B_0 s_j 
				+ \sum_{i=0}^{j-1} \alpha_{3i}(u_{3i}^T s_j) u_{3i}
				+ \alpha_{3i+1}(u_{3i+1}^T s_j) u_{3i+1}
				+ \alpha_{3i+2}(u_{3i+2}^Ts_j)u_{3i+2}.
\end{eqnarray*}
All the terms in the above summation only involve vectors that have been
previously computed.  Having computed $u_{3j+2}$, it is then possible to compute
$\alpha_{3j+1}, \alpha_{3j+2},$ and $u_{3j+1}$.  (The other terms $\alpha_{3j}$
and $u_{3j}$ do not depend on $B_js_j$.)  The following algorithm
computes the terms $\alpha$ and $u$ that are used in Algorithm 4.

\begin{algorithm}[H]
\caption{Unrolling the limited-memory Broyden convex class formula}\label{alg-unrolling}
  \tab \FOR\ $j=0,\ldots, k$
  \tabb $u_{3j} \leftarrow y_j$;
  \tabb $\alpha_{3j}  \leftarrow 1/(s_j^Ty_j)$;
  \tabb $u_{3j+2} \leftarrow B_0s_j + \sum_{i=0}^{j-1} 
  		\left [ \alpha_{3i}(u_{3i}^T s_j) u_{3i}
				+ \alpha_{3i+1}(u_{3i+1}^T s_j) u_{3i+1}
				+ \alpha_{3i+2}(u_{3i+2}^Ts_j)u_{3i+2} \right ]$;
  \tabb $\alpha_{3j+2} \leftarrow -1/(s_j^Tu_{3j+2})$;	
  \tabb $u_{3j+1} \leftarrow \alpha_{3j}y_j + \alpha_{3j+2} u_{3j+2}$;
  \tabb $\alpha_{3j+1} = -\phi/\alpha_{3j+2}$;
  \tab \END\
\end{algorithm}

Assuming $s_j^Ty_j$ is precomputed at each step, Algorithm 5 requires 
$$
	\sum_{j=0}^k 
	\big \{
	3j
	+ 1
	\big \}
	=
	\frac{3(k+1)k}{2}+ k + 1
	= \frac{(k+1)(3k+2)
	}
	{2}
$$
vector inner products and
$$
	\sum_{j=0}^k
	\left \{
	1 + 
	2n + 
	(5n + 3)j
	+ 1
	+ 3n
	+ 1
	\right \} = 
	\frac{(5n+3)(k+2)(k+1)}{2}
$$
additional flops.  
With the $\alpha$'s and the $u$'s computed in Algorithm 5, 
the rest of Algorithm 4 requires 
$$
	\sum_{j=0}^{3(k+1)-1}
	\left \{
		j  + 1 + 1 
	\right \}
	= \frac{3(k+1)(3(k+1)+5)}{2}
$$ 
vector inner products and 
\begin{eqnarray*}
	\sum_{j=0}^{3(k+1)-1}
	\left \{
		n + (2n+1)l
		+ 3
		+ (2n+1)
	\right \}
	&=& \frac{3(k+1)(3(k+1)+1)(2n+1)}{2}  \\
	&& + 3(k+1)(3n+4)
\end{eqnarray*}
additional flops.
Thus, the total number of vector inner products for Algorithm 4 is
$6(k+1)^2 + 7(k+1)$ 
and the total flop count is
$$
	\frac12(k+1)(23kn+52n+12k+42)+ (2n-1)(6(k+1)^2 + 7(k+1)),
$$ which includes
$6(k+1)^2 + 7(k+1)$ vector inner products.

\bigskip

\section{Computing products with $H_{k+1}$ recursively}\label{app:Hk+1}

The inverse of $B_{k+1}$ in \eqref{eqn-broyden} is given by
\begin{equation*}
	H_{k+1}  = H_k  + \frac{1}{s_k^Ty_k}s_ks_k^T 
				- \frac{1}{y_k^TH_k y_k}H_k y_ky_k^TH_k 
				+ \Phi_k (y_k^TH_k y_k)v_kv_k^T,
\end{equation*}
where
\begin{equation*}
	v_k = \frac{s_k}{y_k^Ts_k} - \frac{H_k y_k}{y_k^TH_k y_k}
	\quad
	\text{and} 
	\quad 
	\Phi_k 
	= 
	\frac{(1-\phi)(y_k^Ts_k)^2}{(1-\phi)(y_k^Ts_k)^2 + \phi(y_k^TH_k y_k)(s_k^TB_k s_k)},
\end{equation*}
and $H_k\defined B_k^{-1}$.  Thus, the solution to $B_{k+1}r = z$ can be obtained from 
\begin{eqnarray*}
	B_{k+1}^{-1}z 
	&=&
	\left ( H_0z + \sum_{i=0}^k \frac{1}{s_i^Ty_i}s_is_i^T - \frac{1}{y_i^TH_iy_i}H_iy_iy_i^TH_i
		+ \Phi_i(y_i^TH_iy_i)v_i^Tv_i \right )z
	\\	
	&=& H_0z + \sum_{i=0}^k \frac{s_i^Tz}{s_i^Ty_i}s_i - \frac{y_k^TH_iz}{y_i^TH_iy_i}H_iy_i
		+ \Phi_i(y_i^TH_iy_i)(v_i^Tz)v_i.
\end{eqnarray*}
In order to avoid storing $H_0,\ldots, H_k$ for matrix-vector products,
we make use of the following variables:
\begin{equation}\label{eqn-alphau2}
\begin{array}{rlllrll}
\alpha_{3i} &=& (s_i^Ty_i)^{-1},& \hspace*{1cm}  u_{3i}&=& s_i, \\
\alpha_{3i+2} &=& -(y_i^TH_iy_i)^{-1} &\hspace*{1cm} u_{3i+2}&=& H_iy_i,\\
\alpha_{3i+1} &=&  \Phi_i(y_i^T H_i  y_i ), & \hspace*{1cm} u_{3i+1} &=&
	\displaystyle \frac{s_i}{y_i^Ts_i} - \frac{H_i  y_i}{y_i^TH_i   s_i}.  \\
\end{array}
\end{equation}
With these definitions, we have that
\begin{equation}\label{eqn-Brec}
	B_{k+1}^{-1}z = H_0z + \sum_{i=0}^{k} \alpha_{3i} (u_{3i}^Tz)u_{3i}
	+ \alpha_{3i} (u_{3i+1}^Tz)u_{3i+1}
	+ \alpha_{3i} (u_{3i+2}^Tz)u_{3i+2},
\end{equation}
and, thus, a recursion relation can be used for matrix-vector products
with $H_0,\ldots,H_k$.
Algorithm 6 details how to solve for $r$
in (\ref{eqn-basic}) using the expression for $B_{k+1}^{-1}$ in (\ref{eqn-Brec})
without explicitly storing $H_0,\ldots, H_k$.

\begin{algorithm}[H]
\caption{Computing $r=H_{k+1} z$}\label{alg-5}
  \tab $r\leftarrow H_0z$;
  \tab \FOR\ $j=0,\ldots,3(k+1)-1$
  \tabb \IF\ $\mod(j,3)=0$ 
  \tabbb 	$i\leftarrow j/3$;
  \   	$u\leftarrow s_{i}$;
  \   	$\alpha\gets (y_{i}^Ts_{i})^{-1}$;
  \tabb \ELSEIF\ $\mod(j,3)=1$
  \tabbb $i\leftarrow (j-1)/3$;
  \   $u\gets H_{i}y_{i}$;
  \   $\alpha\gets -(y_{i}^TH_{i}y_{i})^{-1}$;
  \tabb  \ELSE\ 
  \tabbb $i\leftarrow (j-2)/3$;
  \    $u \gets (s_{i})/(y_{i}^Ts_{i}) - (H_{i}y_{i})/(y_{i}^TH_{i}y_{i})$;
  \tabbb  $\Phi \gets (1-\phi)(y_i^Ts_i)^2/\left (
  (1-\phi)(y_i^Ts_i)^2+\phi(y_i^TH_iy_i)(s_i^TB_is_i)
  \right )$;
  \tabbb   $\alpha\gets \Phi(y_{i}^TH_{i}y_{i})$;  
  \tabb \END\ 
  \tabb $r\leftarrow r + \alpha (u^Tz)u$;
  \tab \END\
\end{algorithm}

The following algorithm computes products
involving the matrices $H_i$ for $0 \le i \le k$ for use in Algorithm 6.
The algorithm avoids storing any matrices, and matrix-vector products are
computed recursively. The derivation of this algorithm is similar to that
for Algorithm 5.

\begin{algorithm}[H]
\caption{Unrolling the limited-memory Broyden convex class formula}\label{alg-6}
  \tab \FOR\ $j=0,\ldots, k$
  \tabb $\hat{u}_{3j} \leftarrow s_j$;
  \tabb $\hat{\alpha}_{3j}  \leftarrow 1/(s_j^Ty_j)$;
  \tabb $\hat{u}_{3j+1} \leftarrow H_0y_j + \sum_{i=0}^{j-1} 
  		\left [ \hat{\alpha}_{3i}(\hat{u}_{3i}^T y_j) \hat{u}_{3i}
				+ \hat{\alpha}_{3i+1}(\hat{u}_{3i+1}^T y_j) \hat{u}_{3i+1}
				+ \hat{\alpha}_{3i+2}(\hat{u}_{3i+2}^Ty_j)\hat{u}_{3i+2} \right ]$;
  \tabb $\hat{\alpha}_{3j+1} \leftarrow -1/(y_j^T\hat{u}_{3j+2})$;	
  \tabb $\hat{u}_{3j+2} \leftarrow \hat{\alpha}_{3j}s_j + \hat{\alpha}_{3j+2} \hat{u}_{3j+2}$;
  \tabb Compute $\alpha_{3j+2}$ using Algorithm 5;
  \tabb $\Phi_j \leftarrow  (1-\phi)/( 1-\phi + \phi\hat{\alpha}_{3j}^2/(\hat{\alpha}_{3j+1}\alpha_{3j+2}))$;
  \tabb $\hat{\alpha}_{3j+2} = -\Phi_j/\hat{\alpha}_{3j+1}$;
  \tab \END\

\end{algorithm}

Assuming $s_j^Ty_j$ is precomputed at each step, Algorithm 7 requires 
$$
	\sum_{j=0}^{k}
	\left \{
		3j  + 1 
	\right \}
	+ \frac{(k+1)(3k+2)}{2}
	= (3k+2)(k+1)
$$
vector inner products and
\begin{eqnarray*}
	&& \sum_{j=0}^{k}
	\left \{
		1 +
		(5n+3)j + 
		2n +
		1 +
		3n  + 
		8 +
		1
	\right \} +
	\frac{3(k+1)(3(k+1)+5)}{2}\\
	&& \quad = (5n+3)(k+1)k + (10n+14)(k+1)
\end{eqnarray*} 
additional flops.  
With the $\alpha$'s and the $u$'s computed in Algorithm 7, 
the rest of Algorithm 6 requires 
$3(k+1)$ vector inner products and  
$3(k+1)(2n+1)+n$ additional flops.
Thus, the total flop count for Algorithm 6 is
$$(5n+3)(k+1)k + (k+1)(16n+17) + n + (2n-1)(3k+5)(k+1),$$
which includes $(3k+5)(k+1)$ vector inner products.

\bigskip


\section{Relationships between updates}\label{app:updates}

We note that in the case of the \BFGS{} update ($\phi = 0$),
the compact representation of the BFGS matrix is consistent
with its known compact representation derived in \cite{ByrNS94}.
In particular, $\tilde{M}_k$ in (\ref{eqn-mtilde})
simplifies to 
\begin{eqnarray*}
	\tilde{M}_k 
 & = & (-M_k^{-1} - \Psi_k^TH_0\Psi_k)^{-1}  \\
	&=&
	\begin{bmatrix}
	0 & -R_k - D_k  \\
	-R_k^T - D_k  & -D_k - Y_k^TH_0Y_k 
	\end{bmatrix}^{-1}\\
	&=&
	\begin{bmatrix}
	\bar{R}_k^{-T}(D_k+Y_k^TH_0Y_k)\bar{R}_k^{-1} &  -\bar{R}_k^{-T} \\
	-\bar{R}_k^{-1} & 0 
	\end{bmatrix},
\end{eqnarray*}
where $\bar{R}_k = R_k + D_k$.  Thus, the compact representation
for the inverse of a \BFGS{} matrix is given by
\begin{equation}\label{eqn-bfgsinv}
H_{k+1} =
	H_0 + 
	[ \ S_k \ \ \ H_0 Y_k \ ]
	\begin{bmatrix}
	\bar{R}_k^{-T}(D_k+Y_k^TH_0Y_k)\bar{R}_k^{-1} &  -\bar{R}_k^{-T} \\
	-\bar{R}_k^{-1} & 0 
	\end{bmatrix}
	\begin{bmatrix}
		S_k^T \\
		Y_k^TH_0
	\end{bmatrix},
\end{equation}
which is equivalent to that found in~\cite[Equation (2.6)]{ByrNS94}.

When $\phi = 1$, then $B_{k+1}$ in \eqref{eqn-broyden} is known as the
Davidon-Fletcher-Powell (\DFP{}) update, which preceded the \BFGS{} update.
The \BFGS{} and \DFP{} formulas are known to be \emph{duals} of each other,
meaning one update can be obtained from the other by interchanging $s_k$
with $y_k$ and $B_k$ with $H_k$.  Here, we demonstrate explicitly that the
compact representation of the \BFGS{} and \DFP{} updates are also duals of
each other.  In addition, we show that the compact representation of the
SR1 matrix is self-dual.
%

Consider the compact formulation for the inverse of a \BFGS{} matrix (i.e.,
$\phi=0$) in (\ref{eqn-bfgsinv}).  This is equivalent to
\begin{eqnarray*}
H_{k+1} &=&
	H_0 + 
	[ \ H_0 Y_k \ \ \ S_k \ ]
	\begin{bmatrix}
		0 & -\bar{R}_k^{-1} \\
		 -\bar{R}_k^{-T} & \bar{R}_k^{-T}(D_k+Y_k^TH_0Y_k)\bar{R}_k^{-1} 
	\end{bmatrix}
	\begin{bmatrix}
		Y_k^TH_0\\
		S_k^T 
	\end{bmatrix} \\
	&=& 
	H_0 -
	[ \ H_0 Y_k \ \ \ S_k \ ]
	\begin{bmatrix}
		Y_k^TH_0Y_k + D_k & \bar{R}_k^T \\
		\bar{R}_k & 0
	\end{bmatrix}^{-1}
	\begin{bmatrix}
		Y_k^TH_0\\
		S_k^T 
	\end{bmatrix}. 
\end{eqnarray*}
Now we consider interchanging $B_0$ with $H_0$ and 
$Y_k$ with $S_k$.  Notice that if
$S_k$ and $Y_k$ are interchanged then the upper triangular
part of $S_k^TY_k$ corresponds to the lower triangular part of $Y_k^TS_k$,
implying that $L_k$ and $R_k$ must also be interchanged.
Putting this all together yields:
\begin{eqnarray*}
B_{k+1} &=&
	B_0 -
	[ \ B_0 S_k \ \ \ Y_k \ ]
	\begin{bmatrix}
		S_k^TB_0S_k + D_k & (L_k+D_k)^T \\
		L_k+D_k & 0
	\end{bmatrix}^{-1}
	\begin{bmatrix}
		S_k^TB_0\\
		Y_k^T 
	\end{bmatrix}, 
\end{eqnarray*}
which is the compact formulation for a \DFP{} matrix
(see~\cite[Theorem 1]{ErwJM13}).  In other words, the
compact representations of 
\BFGS{} and \DFP{} are complementary updates of each other.  

For the DFP update, $\Lambda_k = -D_k$  in \eqref{eqn-Lambda} since $\phi = 1$.
Thus, the compact formulation for the inverse of a DFP matrix is given by
\begin{eqnarray*}
	H_{k+1} &=&
	H_0 +
	[ \  S_k \  \ \ H_0Y_k]	
	\begin{bmatrix}
		\ \ D_k & -R_k \\
		-R_k^T & - Y_k^TH_0Y_k		
	\end{bmatrix}^{-1}
	\begin{bmatrix}
		S_k^T \\
		Y_k^TH_0
	\end{bmatrix}\\
	&=&
	H_0  - 
	[ \  H_0Y_k \ \ \ S_k]		
	\begin{bmatrix}
		 Y_k^TH_0Y_k  & \ \ R_k^T  \ \\
		R_k & 	- D_k \  	
	\end{bmatrix}^{-1}
	\begin{bmatrix}
		Y_k^TH_0\\
		S_k^T
	\end{bmatrix}.
\end{eqnarray*}
Interchanging $H_0$ with $B_0$ and $S_k$ with $Y_k$ (and $R_k$ with $L_k$) yields 
$$
	B_{k+1} = 
	B_0 - 
	[ \ B_0S_k \ \ \ Y_k ]
	\begin{bmatrix}
		S_k^TB_0S_k & \ \ L_k^T \ \\
		L_k & -D_k \
	\end{bmatrix}
	\begin{bmatrix}
		S_k^TB_0 \\
		Y_k^T
	\end{bmatrix},
$$
which is the compact formulation of the BFGS matrix (see Eq.\ (2.17) in 
\cite{ByrNS94}.


\medskip

Finally, it is worth noting that the inverse of compact formulation
for an \SR{} matrix is self-dual in the sense of Section 5.2.
That is, replacing $H_0$ with $B_0$ and $S_k$ with $Y_k$ (and, thus,
$R_k$ with $L_k$) in (\ref{eqn-H3})
yields the compact formulation for the \SR{} matrix given by
(\ref{eqn-alt-form-SR1}).


\section{Solves with SR1 matrices.}\label{app:SR1}

The \SR{} update is remarkable in that it
is self-dual: initializing with $B_0^{-1}$ instead of $B_0$, replacing
$s_k$ with $y_k$, and replacing $y_k$ with $s_k$ for all $k$ in
(\ref{eqn-SR1}) results in $B_{k+1}^{-1}$ (see e.g.,~\cite{NocW06}).  Thus,
$r=B_{k+1}^{-1}z$ can be computed using recursion (see Algorithm E.1).  The
first loop of Algorithm E.1 computes $p_i \defined s_i-H_iy_i$ for
$i=0,\ldots k$; the final line of Algorithm E.1 computes
$r=B_{k+1}^{-1}z=B_0^{-1}z + \sum_{i=0}^{k}(p_i^Tz)/(p_i^Ty_i)p_i$.

\begin{algorithm}[H]
\caption{Computing $r=B_{k+1}^{-1} z,$ when $B_{k+1}$
is an \SR{} matrix}\label{alg-sr1-Hz}
\tab \FOR\ $i=0,\ldots k$
\tabb $p_i = s_i-H_0y_i$;
\tabb \FOR $j=0\ldots i-1$
\tabbb $p_i=p_i-((p_j^Ty_i)/(p_j^Ty_j))p_j$;
\tabb \END\
\tab \END
\tab $r\gets H_0z + \sum_{i=0}^k ((p_i^Tz)/(p_i^Ty_i))p_i$;
\end{algorithm}

\bigskip

Thus, in the case of \SR{} updates, solving
linear systems with \SR{} matrices can be performed using vector
inner products; the cost for solving linear systems with an \SR{} matrix is
the same as the cost of computing products with an \SR{} matrix,
which can be done with 
$$
	\left \{
	\sum_{i=0}^k i 
	\right \} +
	k+1
	=
	\frac{(k+2)(k+1)}{2}
$$
vector inner products and 
$$
	\left \{
	\sum_{i=0}^{k}
		2n +
		(2n+1)i
	\right \}
	+ 2n + (k+1)(1 +n)
	= \frac{(k+1)}{2}
	\left (
	k(2n+1)+2(3n+1)
	\right )
	+2n
$$
additional flops. 
It should be noted that unlike members of the restricted Broyden
class, \SR{} matrices can be indefinite, and in particular, numerically
singular.  Methods found in~\cite{Burdakov13,ErwM15} can be used to compute the
eigenvalues (and thus, the condition number) of \SR{} matrices before
performing linear solves to help avoid solving ill-conditioned systems.


\newpage


\section{Calculations in Lemma 1.}\label{app:Lemma1}

The determinant of (\ref{eqn-M0}) is given by
\begin{eqnarray*}
	\tilde{\alpha}_0\tilde{\delta}_0 - \tilde{\beta}_0^2 
	&=&
	-\frac{1-\Phi_0}{(s_0^Ty_0)(y_0^TH_0y_0)}
	-\frac{\Phi_0(1-\Phi_0)}{(s_0^Ty_0)^2}
	-\frac{\Phi_0^2}{(s_0^Ty_0)^2}
	\nonumber
	\\
	&=&
	-\frac{\phi s_0^TB_0s_0}{\Delta_0 s_0^Ty_0}
	-\frac{1-\phi}{\Delta_0}
	\nonumber \\
	&=& \frac{1}{\Delta_0} \left (
	\frac{s_0^TB_0s_0}{\lambda_0}
	\right ), 
\end{eqnarray*}
where $\lambda_0$ in \eqref{eqn-determtilde} is defined in \eqref{eqn-Lambda}.
Thus, the first entry of (\ref{eqn-base0inv}) is given by
\begin{equation*}
	\frac{\tilde{\delta}_0}{\tilde{\alpha}_0\tilde{\delta}_0 - \tilde{\beta}_0^2} 
	=
	-\frac{\phi (s_0^TB_0s_0)}{\Delta_0}
	\frac{\Delta_0\lambda_0}{s_0^TB_0s_0} =
-\phi \lambda_0.
\end{equation*}
The off-diagonal elements of the right-hand side of (\ref{eqn-base0inv})
simplify as follows:
\begin{eqnarray*}
	\nonumber
	\frac{\tilde{\beta}_0}{\tilde{\alpha}_0\tilde{\delta}_0 - \tilde{\beta}_0^2}
	&=&
	- \left (\frac{(1-\phi)(s_0^Ty_0)}{\Delta_0} \right )
	\left (
	{\Delta_0} 
	\left (
	\frac{\lambda_0}{s_0^TB_0s_0}
	\right )
	\right )
	\\ \nonumber
	&=&
	-\frac{(1-\phi)s_0^Ty_0}{-(1-\phi) - \displaystyle \phi \frac{s_0^TB_0s_0}{s_0^Ty_0}}
	\\ \nonumber
	&=&
	-\frac{(1-\phi)(s_0^Ty_0)^2 + \phi s_0^TB_0s_0 (s_0^Ty_0)- \phi s_0^TB_0s_0(s_0^Ty_0)}{-(1-\phi)s_0^Ty_0 - \phi s_0^TB_0s_0}
	\\ \nonumber
	&=&
	s_0^Ty_0  + \phi \frac{s_0^TB_0s_0(s_0^Ty_0)}{-(1-\phi)s_0^Ty_0 - \phi s_0^TB_0s_0}
	\\ 
	&=& 
	s_0^Ty_0 + \phi \lambda_0. 
\end{eqnarray*}
Finally, the last entry of (\ref{eqn-base0inv}) can be simplified as
follows:
\begin{eqnarray*}
	\nonumber
	\frac{\tilde{\alpha}_0}{\tilde{\alpha}_0\tilde{\delta}_0 - \tilde{\beta}_0^2} 
	&=& 
	\frac{1}{s_0^Ty_0} \frac{\Delta_0 \lambda_0}{s_0^TB_0s_0} +
	(1-\phi)y_0^TH_0y_0\frac{\lambda_0}{s_0^TB_0s_0}
	\\ \nonumber
	&=&\frac{\lambda_0}{(s_0^Ty_0)(s_0^TB_0s_0)} \bigg ( \Delta_0 + (1-\phi)(y_0^TH_0y_0) (s_0^Ty_0) \bigg )
	\\ \nonumber
	&=&\frac{1}{-(1-\phi)s_0^Ty_0 - \phi s_0^TB_0s_0} \bigg ( \Delta_0 + (1-\phi)(y_0^TH_0y_0) (s_0^Ty_0) \bigg )
	\\ \nonumber
	&=&
	\frac{(1-\phi)(s_0^Ty_0)^2 + \phi (y_0^TH_0y_0)(s_0^TB_0s_0) +(1-\phi)(y_0^TH_0y_0)(s_0^Ty_0)}
	{-(1-\phi)s_0^Ty_0 - \phi (s_0^TB_0s_0)}
	\\ \nonumber
	&=&
	\frac{(1-\phi)(s_0^Ty_0)^2 + \phi (s_0^TB_0s_0)(s_0^Ty_0) }
	{-(1-\phi)s_0^Ty_0 - \phi (s_0^TB_0s_0)}
	-
	\frac{\phi (s_0^TB_0s_0)(s_0^Ty_0) }
	{-(1-\phi)s_0^Ty_0 - \phi (s_0^TB_0s_0)}
	\\ \nonumber
	& & \quad
	+
	\frac{\phi (y_0^TH_0y_0)(s_0^TB_0s_0) +(1-\phi)(y_0^TH_0y_0)s_0^Ty_0}
	{-(1-\phi)s_0^Ty_0 - \phi (s_0^TB_0s_0)}
	\\ 
	&=& -s_0^Ty_0 - \phi\lambda_0 - y_0^TH_0y_0 . 
\end{eqnarray*}

\bibliographystyle{siam}
\bibliography{CompInv}

\end{document}